\documentclass[11pt,amsfonts]{article}
\usepackage{latexsym}
\usepackage{amssymb}
\usepackage{amsmath}
\usepackage{enumerate}
\usepackage{here} 

\usepackage{layout}
\usepackage{eufrak}

\usepackage{amsthm}
\usepackage[bbgreekl]{mathbbol} 
\usepackage{bm}
\usepackage{bbm}
\usepackage{ulem}  
\usepackage{comment}
\usepackage[dvipdfmx]{graphicx}
\usepackage{latexsym}
\usepackage{amssymb}
\usepackage{amsmath}
\usepackage{enumerate}
\usepackage{color}
\usepackage{layout}
\usepackage{dsfont}
\usepackage{eufrak}
\usepackage{authblk}
\usepackage{here}

\newtheorem{theorem}{Theorem}[section]
\newtheorem{proposition}[theorem]{Proposition}
\newtheorem{lemma}[theorem]{Lemma}
\newtheorem{definition}[theorem]{Definition}

\newtheorem{remark}[theorem]{Remark}

\makeatletter
  
  \@addtoreset{equation}{section}
\makeatother

\def\real{{\mathord{{\rm I\kern-2.8pt R}}}}        
\def\inte{{\mathord{{\rm I\kern-2.8pt N}}}}

\def\sZZ{{\rm Z\kern-2.8ptem{}Z}}

\def\z{{\mathchoice
		{\sZZ}
		{\sZZ}
		{\rm Z\kern-0.30em{}Z}
		{\rm Z\kern-0.25em{}Z} }}
\def\sQQ{{\kern 0.27em \vrule height1.45ex width0.03em depth0em
		\kern-0.30em \rm Q}}
\def\qu{{\mathchoice
		{\sQQ}
		{\sQQ}
		{\kern 0.225em \vrule height1.05ex width0.025em depth0em \kern-0.25em \rm Q}
		{\kern 0.180em \vrule height0.78ex width0.020em depth0em \kern-0.20em \rm Q}
}}
\def\sCC{{\kern 0.27em \vrule height1.45ex width0.03em depth0em
		\kern-0.30em \rm C}}
\def\complex{{\mathchoice
		{\sCC}
		{\sCC}
		{\kern 0.225em \vrule height1.05ex width0.025em depth0em \kern-0.25em \rm C}
		{\kern 0.180em \vrule height0.78ex width0.020em depth0em \kern-0.20em \rm C}
}}

\newcommand{\ba}{\begin{array}}
	\newcommand{\ea}{\end{array}}
\newcommand{\be}{\begin{equation}}
	\newcommand{\ee}{\end{equation}}
\newcommand{\bea}{\begin{eqnarray}}
	\newcommand{\eea}{\end{eqnarray}}
\newcommand{\beaa}{\begin{eqnarray*}}
	\newcommand{\eeaa}{\end{eqnarray*}}

\def\z{\zeta}

\def\l{\lambda}

\font\tenmath=msbm10 \font\sevenmath=msbm7 \font\fivemath=msbm5
\newfam\mathfam \textfont\mathfam=\tenmath
\scriptfont\mathfam=\sevenmath \scriptscriptfont\mathfam=\fivemath

\def \p{I\!\!P}

\def \={{\buildrel {\rm (law)} \over =}}

\def \l{\ell}

\def\D{\Delta}

\newcommand{\basa}{\begin{assumption}}
	\newcommand{\easa}{\end{assumption}}

\newcommand{\bas}{\begin{assum}}
	\newcommand{\eas}{\end{assum}}

\def\bbd{\mathbb d}
\def\limsup{\mathop{\overline{\rm lim}}}

\def\wt{\widetilde}

\def\wt{\widetilde}

\def\wh{\widehat}
\def\wt{\widetilde}
\def\sfB{{\tt B}}
\def\sfF{{\sf F}}
\def\sfa{{\sf a}}
\def\sff{{\sf f}}
\def\sfh{{\sf h}}
\def\sfk{{\sf k}}
\def\sfg{{\sf g}}
\def\sfq{{\sf p}}
\def\sfq{{\sf q}}
\def\ol{\overline}

\def\limsup{\varlimsup}

\newcommand{\ignore}[1]{}
\def\infm{{\infty\text{--}}}
\def\inftym{\infm}
\def\koko{{\coloroy{koko}}}
\def\bd{\begin{description}}
\def\ed{\end{description}}

\def\D2{\bbD_{2,\infty-}}

\def\ba{\bar{A}}

\def\D{{\bf D}}

\def\calc{{\cal C}}
\def\cald{{\cal D}}
\def\cale{{\cal E}}

\def\calh{{\cal H}}

\def\mfh{{\mathfrak H}}

\def\yeq{\>=\>}
\def\yleq{\>\leq\>}
\def\ygeq{\>\geq\>}

\def\sfk{{\sf k}}

\def\sfg{{\sf g}}

\def\sfr{{\sf r}}

\def\simleq{\ \raisebox{-.7ex}{$\stackrel{{\textstyle <}}{\sim}$}\ }

\def\ep{\epsilon}
\def\half{\frac{1}{2}}

\def\y{\vspace*{3mm}\\}
\def\halflineskip{\vspace*{3mm}}
\def\nn{\nonumber}
\def\be{\begin{equation}}
\def\ee{\end{equation}}
\def\bea{\begin{eqnarray}}
\def\eea{\end{eqnarray}}
\def\beas{\begin{eqnarray*}}
\def\eeas{\end{eqnarray*}}
\def\bi{\begin{itemize}}
\def\ei{\end{itemize}}
\def\im{\item}
\def\bd{\begin{description}}
\def\ed{\end{description}}
\def\l{\left}
\def\r{\right}

\newcommand{\bbB}{{\mathbb B}}

\newcommand{\bbD}{{\mathbb D}}

\newcommand{\bbF}{{\mathbb F}}

\newcommand{\bbI}{{\mathbb I}}

\newcommand{\bbN}{{\mathbb N}}

\newcommand{\bbR}{{\mathbb R}}
\newcommand{\bbS}{{\mathbb S}}

\newcommand{\bbW}{{\mathbb W}}

\newcommand{\bbY}{{\mathbb Y}}
\newcommand{\bbZ}{{\mathbb Z}}

\newif\ifcol
\colfalse 
\ifcol

\newcommand{\colorr}{\color[rgb]{0.8,0,0}}
\newcommand{\colorg}{\color[rgb]{0,0.5,0}}

\else
\newcommand{\colorr}{\color{black}}
\newcommand{\colorg}{\color{black}}
\fi
\newif\ifcol
\colfalse 
\ifcol

\newcommand{\sred}{\color[rgb]{0.8,0,0}}
\newcommand{\sblue}{\color[rgb]{0,0,0.8}}
\else
\newcommand{\sred}{\color{black}}
\newcommand{\sblue}{\color{black}}
\fi
\newif\ifcol
\colfalse 
\ifcol

\else
\fi
\newif\ifcol
\colfalse 
\ifcol
\newcommand{\vred}{\color[rgb]{0.8,0,0}}

\renewcommand{\koko}{{\color{yellow} kokoko}}

\else
\newcommand{\vred}{\color{black}}
\renewcommand{\koko}{{\color{yellow} kokoko}}
\fi
\newif\ifcol
\colfalse 
\ifcol
\newcommand{\xred}{\color[rgb]{0.8,0,0}}

\renewcommand{\koko}{{\color{yellow} kokoko}}

\else
\newcommand{\xred}{\color{black}}
\renewcommand{\koko}{{\color{yellow} kokoko}}
\fi
\newif\ifcol
\colfalse 
\ifcol
\newcommand{\ared}{\color[rgb]{0.8,0,0}}

\newcommand{\ablue}{\color[rgb]{0,0,0.8}}
\renewcommand{\koko}{{\color{yellow} kokoko}}

\else
\newcommand{\ared}{\color{black}}
\newcommand{\ablue}{\color{black}}
\renewcommand{\koko}{{\color{yellow} kokoko}}
\fi
\newif\ifcol
\colfalse 
\ifcol
\newcommand{\bred}{\color[rgb]{0.8,0,0}}

\renewcommand{\koko}{{\color{yellow} kokoko}}

\else
\newcommand{\bred}{\color{black}}
\renewcommand{\koko}{{\color{yellow} kokoko}}
\fi
\excludecomment{en-text}
\includecomment{jp-text}
\excludecomment{comment}
\newcommand{\gray}{\color[rgb]{0.5,0.5,0.5}}

\begin{document}

\renewcommand{\thefootnote}{\fnsymbol{footnote}}

\renewcommand{\thefootnote}{\fnsymbol{footnote}}

\title{
Asymptotic  expansion of the drift estimator for the fractional Ornstein-Uhlenbeck process 
\footnote{
This work was in part supported by 
Japan Science and Technology Agency CREST JPMJCR14D7, JPMJCR2115;  
Japan Society for the Promotion of Science Grants-in-Aid for Scientific Research 
Nos. 17H01702, 23H03354 (Scientific Research);  
and by a Cooperative Research Program of the Institute of Statistical Mathematics. 
C. Tudor also acknowledges partial support from Labex CEMPI (ANR-11-LABX-007-01), ECOS SUD (project C2107), and from the Ministry of Research, Innovation and Digitalization (Romania), grant CF-194-PNRR-III-C9-2023.
}
}
\author[1]{Ciprian A. Tudor}
\author[2,3,4]{Nakahiro Yoshida}
\affil[1]{Universit\'e de Lille 1
\footnote{Universit\'e de Lille 1: 59655 Villeneuve d'Ascq, France}
        }
\affil[2]{Graduate School of Mathematical Sciences, University of Tokyo
\footnote{Graduate School of Mathematical Sciences, University of Tokyo: 3-8-1 Komaba, Meguro-ku, Tokyo 153-8914, Japan. e-mail: nakahiro@ms.u-tokyo.ac.jp}
        }
\affil[3]{CREST, Japan Science and Technology Agency}
\affil[4]{The Institute of Statistical Mathematics
        }
\maketitle

\begin{abstract}\noindent
{\sblue
We present an asymptotic expansion formula of an estimator for the drift coefficient of the fractional Ornstein-Uhlenbeck process. As the machinery, we apply the general expansion scheme for Wiener functionals recently developed by the authors  \cite{tudor2023high}. 
The central limit theorem in the principal part of the expansion has the classical scaling $T^{1/2}$. 
However, the asymptotic expansion formula is a complex 
in that 
the order of the correction term becomes the classical $T^{-1/2}$ for $H\in(1/2,5/8)$, but 
$T^{4H-3}$ for $H\in[5/8,3/4)$. 
}
\end{abstract}

\vskip0.3cm
\noindent
{\bf 2010 AMS Classification Numbers: } 62M09, 60F05, 62H12


\section{Asymptotic expansion of an estimator for a fractional Ornstein-Uhlenbeck process }
We consider the Langevin equation 
\begin{equation}\label{1}
\left\{
\begin{array}{ccl}
		dX_{t}&=&-\theta X_tdt + \sigma dB_{t}, \hskip0.5cm t\geq 0,\\
		X_0&=&x_0,
\end{array}\right.
\end{equation}
where $x_0$ is a constant and 
$\left( B_{t}, t\geq 0\right) $ is a fractional Brownian motion with Hurst index $H\in (1/2,1)$. 
Suppose that the parameter space $\Theta$ is a bounded open set in $\bbR$ 
satisfying $\ol{\Theta}\subset(0,\infty)$, and that 
the true value of $\theta$ is in $\Theta$. 
In what follows, the true value is also denoted by $\theta$ for notational simplicity.

From (\ref{1}), 
\bea\label{202308071211}
X_t &=& e^{-\theta t}x_0+\int_0^te^{-\theta(t-s)}\sigma dB_s, 
\eea
where the stochastic integral is regarded as a 
{\sblue Wiener integral, 
i.e., an divergence integral with respect to the fractional Brownian motion $B$.}

Hu and Nualart \cite{hu2010parameter} investigated the estimator $\wt{\theta}_T$ defined by 
\bea\label{202308062338}
\wt{\theta}_T
&=&
\bigg(\frac{1}{\sigma^2H\Gamma(2H)T}\int_0^TX_t^2\bigg)^{-\frac{1}{2H}}. 
\eea
In the inferential theory, the estimator $\wt{\theta}_T$ is regarded as an M-estimator for the estimating equation 
\bea\label{202308071656}
\int_0^TX_t^2dt-\wt{\nu}_T(\vartheta)\yeq0
\eea
for
\bea
\wt{\nu}_T(\vartheta)
&=& 
\mu(\vartheta)T\quad\text{with}\quad
\mu(\vartheta)\yeq
\sigma^2H\Gamma(2H)\vartheta^{-2H}. 
\eea

We remark that $\wt{\theta}_T$ is an approximately moment estimator but not the exact moment estimator since 
$\ol{\nu}_T(\theta):=E\big[\int_0^TX_t^2dt\big]$ is decomposed as 
$\ol{\nu}_T(\theta)=\wt{\nu}_T(\theta)+\ol{b}_T(\theta)$ and $\ol{b}_T(\theta)$ does not vanish 
though it is of {\sblue order of} $O(1)$ as $T\to\infty$, according to Lemma \ref{202308132345}. 
Since it is common to use a bias-corrected estimator in the higher-order inference, 
we will consider the estimator 
\beas
\wh{\theta}_T^o
&=&
\wt{\theta}_T-{\sred T^{-\half-\sfq(H)}}\beta\big(\wt{\theta}_T\big), 
\eeas
where $\beta=\beta_H\in C^\infty_B(\Theta)$, i.e., $\beta$ is smooth on $\Theta$ and all its derivatives are bounded on $\Theta$, {\sred and $\sfq=\sfq(H)$ is a number define by (\ref{202307110532}).} 
The value of $\wh{\theta}_T^o$ can exceed the boundary of $\Theta$, not necessarily due to the $\beta$ term, therefore 
the estimator $\wh{\theta}_T$ we will consider is more precisely defined as 
\bea\label{202308132357}
\wh{\theta}_T
&=& 
\left\{
\begin{array}{cl}
\wh{\theta}_T^o&\text{ if }\wt{\theta}_T\in\Theta\text{ and }\wh{\theta}_T^o\in\Theta,\y
\theta_*&\text{ otherwise},
\end{array}\right.
\eea
where $\theta_*$ is a prescribed value in $\Theta$. 
The choice of the value $\theta_*$ will not affect asymptotically in any order of expansion.

Hu and Nualart \cite{hu2010parameter} proved that, for $H\in \left( \frac{1}{2}, \frac{3}{4}\right)$,
\begin{equation*}
		\sqrt{T} \big( \wt{\theta}_{T}-\theta \big) \to^d N(0,c_0)
\end{equation*}
as $T\to\infty$, with $c_0$ defined as (\ref{202308141446}). 
{\ared On the other hand, Hu et al. showed in [5] that the estimator (1.6) converges to a non-Gaussian distribution (a Rosenblatt random variable), when $H\in (3/4, 1)$.  Other estimators for the drift parameter of the fractional Ornstein-Uhlenbeck process  have  been analyzed, among others, in 
Brouste and Kleptsyna \cite{brouste2010asymptotic}, 
Chen and Li \cite{chen2021berry}, 
Cheng and Zhou \cite{chen2021parameter}, and 
El Onsy, Es-Sebaiy and Viens \cite{el2017parameter}. 
}

In this paper, we will give an asymptotic expansion for the distribution of $\sqrt{T} \big( \widehat{\theta}_{T}-\theta \big) $. 
The order  $\sfq$ of the expansion is defined as 
\bea\label{202307110532}
\sfq\yeq\sfq(H)
&=&
\l\{\begin{array}{ll}
\half&\big(H\in\big(\frac{1}{2},\frac{5}{8}\big]\big)\y
-4H+3&\big(H\in\big(\frac{5}{8},\frac{3}{4}\big)\big)
\end{array}\r.
\eea

{\sblue 
The $k$-th Hermite polynomial $H_k(x;0,c_0)$ is defined by 
\beas 
H_k(x;0,c_0)
&=& 
e^{2^{-1}c_0^{-1}x^2}(-\partial_x)^ke^{-2^{-1}c_0^{-1}x^2}
\qquad(x\in\bbR).
\eeas
}
We consider the approximate density function 
\bea\label{202308141943}
p_{H,T}(x)
&=& 
\phi(x;0,c_0)\bigg(1
+1_{\{H\in[\frac{5}{8},\frac{3}{4})\}}2^{-1}c_2H_2(x;0,c_0)T^{4H-3}
\nn\\&&\hspace{63pt}
+1_{\{H\in(\half,\frac{5}{8}]\}}3^{-1}c_3H_3(x;0,c_0)T^{-\half}
\nn\\&&\hspace{63pt}
+
c_1H_1(x;0,c_0){\sred T^{-\sfq(H)}}
\bigg), 
\eea
where the constants $c_0,...,c_3$ {\sred depending on $H$ and $\theta$} will be specified later 
at (\ref{202308141446}) and (\ref{202308141938}). 
For $a,b>0$, we denote by $\cale(a,b)$ the set of measurable functions $g:\bbR\to\bbR$ such that 
$|g(x)|\leq a(1+|x|^b)$ for all $x\in\bbR$. 
The main theorem of this paper is here. 
\begin{theorem}\label{202308141824}
Suppose that $H\in(1/2,3/4)$. Then 
\bea\label{202403210322}
\sup_{g\in\cale(a,b)}\bigg|E\big[g\big(T^{1/2}(\wh{\theta}_T-\theta)\big)\big]
-
\int_\bbR g(x)p_{H,T}(x)dx\bigg|
&=&
o(T^{-\sfq(H)})
\eea
as $T\to\infty$, 
for every $a,b>0$. 
\end{theorem}
The function $\beta$ set so as to satisfy $c_1=0$ corrects the second-order bias. 
{\bred 
In Section \ref{202403280137}, the real performance of the formula $p_{H,T}$ will be investigated 
in several cases by simulations. }

{\ablue

We will treat {\bred mainly} the asymptotic expansion formula (\ref{202308141943}) with the threshold $5/8$ changing the shape of the formula 
by the indicator functions. The expansion formula is still valid even if we remove the indicator functions 
and keep all terms 
{\bred because the exponents of $T$ automatically count the order of terms and the smaller terms, even if they remain in the formula, do not affect the error bound for a given value of $H$. 
More precisely, 
\begin{theorem}\label{202403210212}
Suppose that $H\in(1/2,3/4)$. 
Then there exist constants $c_0,c_{1,1}^+,c_{1,2}^+,c_2,c_3$ such that 
for 
\bea\label{202403210310}
p_{H,T}^+(x)
&=& 
\phi(x;0,c_0)\bigg(1
+2^{-1}c_2H_2(x;0,c_0)T^{4H-3}
+3^{-1}c_3H_3(x;0,c_0)T^{-\half}
\nn\\&&\hspace{63pt}
+
c_{1,1}^+H_1(x;0,c_0)T^{-\half}
+c_{1,2}^+H_1(x;0,c_0)T^{-\sfq(H)}
\bigg), 
\eea
it holds that 
\bea\label{202403210214}
\sup_{g\in\cale(a,b)}\bigg|E\big[g\big(T^{1/2}(\wh{\theta}_T-\theta)\big)\big]
-
\int_\bbR g(x)p_{H,T}^+(x)dx\bigg|
&=&
o(T^{-\sfq(H)})
\eea
as $T\to\infty$, 
for every $a,b>0$. 
\end{theorem}

The constants $c_0$, $c_2$, $c_3$ are the same as those of $p_{H,T}$. 
The constants $c_{1,1}^+$ and $c_{1,2}^+$ are given in (\ref{202403210242}). 
}

In asymptotic expansions in general, such a ``redundant'' formula {\bred may give} in practice a better approximation to the distribution though there is no theoretical explanation except for an intuition that such a primitive formula has more information than the slimmed formulas obtained by further neglecting smaller terms. 
}

{\ared Concluding this section, here are several comments.} %
Hu, Nualart and Zhou \cite{hu2019parameter} presented limit theorems for general Hurst parameter. 
The Berry-Esseen bound for the parameter estimation is discussed, {\ared among others,} in 
Kim and Park \cite{kim2017optimal}, 
Chen, Kuang and Li \cite{chen2020berry}, 
{\ared and Chen and Li \cite{chen2021berry}. }

For estimation of the Hurst coefficient, we refer the reader to 
Istas and Lang \cite{istas1997quadratic}, 
Kubilius and Mishura \cite{kubilius2012rate}, 
Kubilius, Mishura and Ralchenko \cite{kubilius2017parameter} and 
Berzin, Latour and Le{\'o}n \cite{berzin2014inference}. 
Asymptotic expansions are discussed in 
Mishura, Yamagishi and Yoshida \cite{mishura2023asymptotic}. 
A related expansion for the quadratic form for a stochastic differential equation driven by a fractional Brownian motion 
(in particular for the estimator for a constant volatility parameter) 
is in Yamagishi and Yoshida \cite{yamagishi2022order, yamagishi2023order}. 
Tudor and Yoshida \cite{tudor2020asymptotic} gave asymptotic expansion of the quadratic variation of a mixed fractional Brownian motion. 

In this article, we consider an asymptotic expansion for a fractional process, while this problem has been studied 
well for diffusion processes: 
Mykland \cite{Mykland1992}, 
Yoshida \cite{Yoshida1997, Yoshida2004}, 
Kusuoka and Yoshida \cite{KusuokaYoshida2000}, 
Sakamoto and Yoshida \cite{SakamotoYoshida2003,SakamotoYoshida2004,SakamotoYoshida2009}
and 
Kutoyants and Yoshida \cite{KutoyantsYoshida2007}, 
just to mention a few. 

The general expansion formula by Tudor and Yoshida \cite{tudor2023high} was applied in this article. 
A different formulation using a limit theorem to specify the correction term 
is in Tudor and Yoshida \cite{tudor2019asymptotic}. 

{\ared 
The following sections are devoted to the proof of Theorem \ref{202308141824}. 
The asymptotic expansion formula is specified with the Gamma factors defined in Section \ref{202402261415}. 
Since the stochastic expansion of the error of the estimator will be expressed with certain basic variables, 
we derive expansions for their Gamma factors in Section \ref{202402261411}. 
From these expansions, Section \ref{202402261422} gives an asymptotic expansion of the sum $\bbS_T$ 
of the basic variables (Proposition \ref{202308141705}). 
In Section \ref{202402261439}, we obtain a stochastic expansion of the error of the estimator 
by using $\bbS_T$ (Equation (\ref{202308141731})), and in Section \ref{202402261518}, it will be used to prove Theorem \ref{202308141824}, 
with the aid of the perturbation method.} 
{\bred Theorem \ref{202403210212} is proved by a minor change of that of Theorem \ref{202308141824}.}

\begin{en-text}
We will consider a moment estimation type of M-estimator defined by the estimating equation 
\bea\label{202308071656}
\psi_T(\theta) 
\yeq
\int_0^TX_t^2dt-{\nu}_T(\theta)\yeq0,
\eea
where 

{\sblue
The estimator $\wh{\theta}_T^{o}$ is the M-estimator characterized by (\ref{202308071656}) 
with $\nu_T$ given by

Another natural choice of the function $\nu_T$ is the expectation of $\int_0^TX_t^2dt$ itself: 

Suppose that 
\bea\label{202308080141}
\nu_T(\theta) &=& \wt{\nu}_T(\theta)
+b_T(\theta)
\eea
with a function $b_T\in C^\infty(\Theta)$ such that for some positive number $T_0$, for all $i\in\bbZ_+$, 
\bea\label{202308080142}
\sup_{\theta\in\Theta,T>T_0} \big|\partial_\theta^ib_T(\theta)\big| &<& \infty
\eea
and that 
\bea\label{202308080143}
\lim_{T\to\infty}\partial_\theta^ib_T(\theta) &\text{exists}=:& b^{(i)}_\infty(\theta)
\eea
for a function $b_\infty:\Theta\to\bbR$. 
Under  (\ref{202308080142}) and (\ref{202308080143}), 
the function $b_\infty$ is continuous and 
the convergence (\ref{202308080143}) is uniform in $\theta\in\Theta$. 

Any measurable map $\wh{\theta}_T:C([0,T])\to\ol{\Theta}$ is called a minimum contrast estimator for $\theta$ 
with respect to the contrast function $|\psi_T|$ if 
\beas 
\big|\psi_T\big(\wh{\theta}_T)\big|
&=&
\inf_{\theta\in\ol{\Theta}}|\psi_T(\theta)|. 
\eeas
We consider an M-estimator $\wh{\theta}_T$ for $\theta$ with respect to $\psi_T$. 
As will be seen, $P\big[\psi_T\big(\wh{\theta}_T\big)=0]\to1$ as $T\to\infty$. 
{\gray 
In other words, $\wh{\theta}_T$ is an M-estimator for the estimating function $\psi_T(\theta)$. }
\end{en-text}

\section{Gamma factors and their representations}\label{202402261415}
{\ared 
To get the asymptotic expansion (\ref{202308141943}) of the estimator $\wh{\theta}_T$ of (\ref{202308132357}), we will use the method developed in Tudor and Yoshida \cite{tudor2023high}, which is based on the analysis of its  gamma factors. Therefore, we introduce below these random variables and then we study their asymptotic behavior in the later sections for the functionals associated with the stochastic expansion of $\wh{\theta}_T$. 
}

To accommodate a fractional Brownian motion, 
prepare the set $\cale$ of step functions on $\bbR_+=[0,\infty)$, and 
introduce an inner product on $\cale$ such that 
\beas 
\langle 1_{[0,t]},1_{[0.s]}\rangle_\calh
\yeq
R_H(t,s)
:=\>
\half\big(t^{2H}+s^{2H}-|t-s|^{2H}\big)
\eeas
for $t,s\in\bbR_+$. 
Define the Hilbert space $\calh$ as the closure of $\cale$ with respect to 
$\|\cdot\|_\calh=\langle\cdot,\cdot\rangle_\calh^{1/2}$. 
In the case $H\in(1/2,1)$, the space $\calh$ has a subspace $|\calh|$ 
of all measurable functions $\sfh:\bbR_+\to\bbR$ satisfying 
\beas 
\int_0^\infty\int_0^\infty|\sfh_t||\sfh_s||t-s|^{2H-2}dsdt &<& \infty.
\eeas
For elements $\sfh,\sfg\in|\calh|$, 
\beas 
\langle\sfh,\sfg\rangle_\calh&=&\alpha_H\int_0^\infty\int_0^\infty\sfh_t\sfg_s|t-s|^{2H-2}dsdt, 
\quad \alpha_H\yeq H(2H-1). 
\eeas

We consider an isonormal Gaussian process $\bbW=\big(\bbW(\sfh)\big)_{\sfh\in\calh}$ on the Hilbert space $\calh$. 
Then, $B_t=\bbW(1_{[0,t]})$ ($t\in\bbR_+$) form a fractional Brownian motion with the Hurst coefficient $H$. 
We will apply the Malliavin calculus associated with $\bbW$. 
{\xred We denote the Malliavin derivative by $D$, and the Malliavin operator by $L$. 
See Nualart \cite{Nualart2006}, Nourdin and Peccati \cite{nourdin2012normal} and Ikeda and Watanabe \cite{IkedaWatanabe1989} 
for the concepts of the Malliavin calculus.} 

\begin{en-text}
$|\calh|$: 
the set of measurable functions $\sfh:\bbR_+\to\bbR$ satisfying 
$\langle\sfh,\sfh\rangle_{\calh}<\infty$, where 
\beas 
\langle \sfh,\sfg\rangle_{\calh} &=& 
\alpha_H\int_{\bbR_+^2}\sfh_t\sfg_s|t-s|^{2H-2}dsdt\quad (\sfh,\sfg\in\calh)
\eeas
for $H>1/2$. 
\end{en-text}
\begin{en-text}
Suppose that $F=I_2(f)$ for $f\in \calh^{\odot2}$. Then 
\beas 
E[\Gamma^{(2)}(F)]
&=& 
2^{-1}E[\|DF\|_\calh^2]
\yeq
2\|f\|_{\calh^{\otimes2}}^2
\eeas
\beas 
E[\Gamma^{(3)}(F)]
&=&
4\langle f\otimes_1f,f\rangle_{\calh^{\otimes2}}
\eeas
\beas 
E[\Gamma^{(4)}(F)]
&=&
8\langle f\otimes_1f\otimes_1f,f\rangle_{\calh^{\otimes2}}
\eeas
\end{en-text}

For $\sfF=(\sfF_i)_{i=1,...,d}\in{\bbD_{1,2}}^d$, the gamma factors $\Gamma^{(m)}(\sfF_{i_1},...,\sfF_{i_m})$ 
for $(i_1,...,i_m)\in\{1,...,d\}^m$ are defined as 
\beas 
\Gamma^{(1)}(\sfF_{i_1})&=&\Gamma^{(1)}_{i_1}(\sfF) \yeq \sfF_i-E[\sfF_{i_1}],\y
\Gamma^{(m)}(\sfF_{i_1},...,\sfF_{i_m}) &=& \Gamma^{(m)}_{i_1,...,i_m}(\sfF) 
\yeq \big\langle D(-L)^{-1}\Gamma^{(m-1)}(\sfF_{i_1},...,\sfF_{i_{m-1}}),D\sfF_{i_m}\big\rangle_\calh\quad(m\geq2). 
\eeas
The map $(\sfF_{i_1},...,\sfF_{i_m})\mapsto\Gamma^{(m)}(\sfF_{i_1},...,\sfF_{i_m})$ is multi-linear. 
Tudor and Yoshida \cite{tudor2023high} used the notation $\Gamma^{(m)}_{i_1,...,i_m}(\sfF)$ 
for $\Gamma^{(m)}(\sfF_{i_1},...,\sfF_{i_m})$. 
The second gamma factor 
$\Gamma^{(2)}(\sfF_{i_1},\sfF_{i_2})$ is in general different from the carr\'e du champ 
$\Gamma(\sfF_{i_1},\sfF_{i_2})=\langle \sfF_{i_1},\sfF_{i_2}\rangle_\calh$.

Suppose that {\sblue a $d$-dimensional random variable $\sfF=(\sfF_i)_{i=1,...,d}$ has} the representation
\bea\label{202307030126}
\sfF_i=I_2(\sff_i)+c_i
\eea
for some $\sff_i\in\calh^{\odot2}$ and $c_i\in\bbR$. 
In this special case, the gamma factors have the following expressions: 
\beas 
\Gamma^{(1)}(\sfF_{i_1})&=& \sfF_{i_1}-c_{i_1},\\
\Gamma^{(2)}(\sfF_{i_1},\sfF_{i_2}) &=& 
2\langle I_1(\sff_{i_1}),I_1(\sff_{i_2})\rangle_\calh
\yeq
2I_2(\sff_{i_1}\otimes_1\sff_{i_2})+2\langle \sff_{i_1},\sff_{i_2}\rangle_{\calh^{\otimes2}}
\nn\\
\Gamma^{(3)}(\sfF_{i_1},\sfF_{i_2},\sfF_{i_3})&=&
2^2\langle I_1(\sff_{i_1}\otimes_1\sff_{i_2}),I_1(\sff_{i_3})\rangle_\calh
\nn\\&=&
2^2I_2\big(\sff_{i_1}\otimes_1\sff_{i_2}\otimes_1\sff_{i_3}\big)
+2^2\big\langle \sff_{i_1}\otimes_1\sff_{i_2},\sff_{i_3}\big\rangle_{\calh^{\otimes2}}.
\eeas
Generally, 
\bea\label{202307030127}
\Gamma^{(m)}(\sfF_{i_1},...,\sfF_{i_m}) &=& 
2^{m-1}I_2\big(\sff_{i_1}\otimes_1\cdots\otimes_1\sff_{i_m}\big)
+2^{m-1}\big\langle \wt{\sff_{i_1}\otimes_1\cdots\otimes_1\sff_{i_{m-1}}},\sff_{i_m}\big\rangle_{\calh^{\otimes2}}
\eea
for $(i_1,...,i_m)\in\{1,...,d\}^m$ and $\sfF_i$ of (\ref{202307030126}), 
where $\wt{\ }$ means the symmetrization. 

\section{Estimates of the gamma factors of the basic variables}\label{202402261411}
\subsection{Basic variables}

Let 
\beas 
u_T(s,t) &=&  
K_UT^{-1/2}e^{-\theta|s-t|}1_{[0,T]^2}(s,t)
\text{ with }K_U \yeq 
-\frac{\theta^{2H}}{4H^2\Gamma(2H)}, 
\nn\\
v_T(s,t) &=& 
K_VT^{-1/2}e^{-\theta(T-s)-\theta(T-t)}1_{[0,T]^2}(s,t)
\text{ with }K_V \yeq \frac{\theta^{2H}}{{\vred4}H^2\Gamma(2H)}, 
\nn\\
w_T(t) &=& 
K_WT^{-1/2}(e^{-\theta t}-e^{-2\theta T+\theta t})1_{[0,T]}(t)
\text{ with }K_W \yeq 
-\frac{x_0\theta^{2H}}{2\sigma H^2\Gamma(2H)}. 
\eeas
We will treat the multiple integrals
\bea\label{202308091408}
U_T=I_2(u_T),\quad V_T=I_2(v_T) \ \text{ and }\ W_T=I_1(w_T). 
\eea 
{\xred 
These variables will play an important role in this article to derive the asymptotic expansion. 
In fact, the estimator $\wh{\theta}_T$ will be related with the sum of them in (\ref{202308131614}). 
}

\subsection{Gamma factors of $U_T$ and $V_T$}
Since $U_T$ and $V_T$ have the form of (\ref{202308091408}), the formula (\ref{202307030127}) gives 
\bea\label{202307030144}
\Gamma^{(m)}(\sfF_T,...,\sfF_T) 
&=& 
2^{m-1}I_2\big(\underbrace{\sff_T\otimes_1\cdots\otimes_1\sff_T}_{m}\big)
+2^{m-1}\big\langle \underbrace{\sff_T\otimes_1\cdots\otimes_1\sff_T}_{m-1},\sff_T\big\rangle_{\calh^{\otimes2}}
\eea
for $m\geq2$ and $\sfF_T=I_2(\sff_T)=U_T$ and $V_T$ with $\sff_T=u_T$ and $v_T$, respectively.

\subsection{Estimates for $E[\Gamma^{(m)}(U_T,...,U_T)]$}

Let 
\bea\label{202307060903}
a(x_1,x_2,x_3) 
&=&
e^{-\theta|x_1-x_2|}|x_2-x_3|^{2H-2}
\eea
for $x_1,x_2,x_3\in\bbR$, and 
\beas
\ol{a}(x) 
\yeq
\int_\bbR e^{-\theta|z|}|z-x|^{2H-2}dz
\eeas
for $x\in\bbR$. 
Then 
\bea\label{202307060939}
\ol{a}(x)=\ol{a}(|x|)\quad\text{and}\quad
\ol{a}(x-y) \yeq
\int_\bbR a(x,z,y)dz
\ygeq
\int_A a(x,z,y)dz
\eea
{\sblue for any $x,y\in\bbR$ and any one-dimensional Borel set $A$.}
The functions $a(x_1,x_2,x_3)$ and $\ol{a}(x)$ depend on $\theta$ and $H$. 
\begin{en-text}
\beas 
\ol{a}(x-y) 
&=& 
\int_\bbR e^{-\theta|x-z|}|z-y|^{2H-2}dz
\yeq
\int_\bbR e^{-\theta|z|}|z-r|^{2H-2}dz
\yeq
\ol{a}(r)\quad(r=|x-y|)
\eeas
for $x,y\in\bbR$. 
\end{en-text}
\begin{lemma}\label{202307060258}
There exists a positive constant $C$ depending on $(\theta,H)$, such that 
\bea\label{202307050606}
\ol{a}(r)\leq C(1\wedge r^{2H-2})\quad (\forall r\geq0). 
\eea
\end{lemma}
\proof 
Notice that $2|z|\geq1$ for $|z-1|\leq1/2$. 
For $r>0$, we have 
\beas
\ol{a}(r) 
&=&
r^{2H-2}\int_\bbR re^{-\theta r|z|}|z-1|^{2H-2}dz
\nn\\&\leq&
r^{2H-2}\bigg(2^{2-2H}\int_{\{z:|z-1|>1/2\}} re^{-\theta r|z|}dz+\int_{\{z:|z-1|\leq1/2\}} 
\sup_{z'\in\bbR}\big(2|z'|\>re^{-\theta r|z'|}\big)|z-1|^{2H-2}dz\bigg)
\nn\\&\leq&
2^{3-2H}\theta^{-1}\big(1+(2H-1)^{-1}e^{-1}\big)r^{2H-2}\quad\text{since }H>1/2, 
\nn
\eeas
besides, 
$
\ol{a}(r) 
\leq
\int_{\{z:|z-r|\geq1\}}e^{-\theta|z|}dz+\int_{\{z:|z-r|<1\}}|z-r|^{2H-2}dz
<
2\theta^{-1}+2(2H-1)^{-1}<\infty
$. 
\qed\halflineskip

Here is a common estimate for a multiple integral. 
\begin{lemma}\label{202308120413}
Let $m\geq2$ and $H\in \left( \frac{1}{2}, \frac{m+1}{2m}\right)$. 
Suppose that functions $\alpha_i:\bbR\to\bbR_+$ $(i=1,...,m)$ 
satisfy
\bea\label{202308120433}
\alpha_i(x)&\leq& C(1\wedge|x|^{2H-2})\quad(x\in\bbR)
\eea
for some positive constant $C$. 
Then 
\beas 
\int_{\bbR^{m-1}}\alpha_1(x_1)\alpha_2(x_1-x_2)\cdots
\alpha_{m-1}(x_{m-2}-x_{m-1})\alpha_m(x_{m-1})dx_1\cdots dx_{m-1}
&<& \infty.
\eeas
\end{lemma}
\proof
By Young's inequality and H\"older's inequality, we obtain 
\bea\label{202308091838}&&
\int_{\bbR^{m-1}}\alpha_1(x_1)\alpha_2(x_1-x_2)\cdots
\alpha_{m-1}(x_{m-2}-x_{m-1})\alpha_m(x_{m-1})dx_1\cdots dx_{m-1}
\nn\\&=&
\big\|(\alpha_1*\cdots*\alpha_{m-1})\times\alpha_m\big\|_{L^1(\bbR)}
\yleq 
\prod_{i=1}^m
\|\alpha_i\|_{L^{\frac{m}{m-1}}}.
\eea
Since $H<\frac{m+1}{2m}$, we have 
$(2H-2) \frac{m}{m-1}<-1$, and hence 
$\|\alpha_i\|_{L^{\frac{m}{m-1}}}<\infty$ from the inequality (\ref{202308120433}). 
\qed\halflineskip

Let 
\beas 
{\vred C_U(m,H,\theta)}
&=&
{\vred 2^mmK_U^{m}\alpha_H^m}\int_{(0,\infty)^{2m-1}}
a(0, x_{2}, x_{3}) a(x_{3}, x_{4}, x_{5})\cdots 
\nn\\&&\hspace{50pt}\cdots
a(x_{2m-3},x_{2m-2}, x_{2m-1})a(x_{2m-1}, x_{2m},0)\>dx_{2}\cdots dx_{2m}. 
\eeas
According to Hu and Nualart \cite{hu2010parameter}, 
\beas
\int_{(0,\infty)^3}a(0,x_2,x_3)a(x_3,x_4,0)dx_2dx_3dx_4
&=&
\theta^{-4H+1}d_H
\eeas
for
\beas 
d_H 
&=& 
(4H-1)\bigg\{\frac{\Gamma(2H-1)^2}{2}+\frac{\Gamma(2H-1)\Gamma(3-4H)\Gamma(4H-2)}{\Gamma(2-2H)}\bigg\}.
\eeas
Therefore,
\bea\label{202308150202}
C_U(2,H,\theta)
&=& 
\frac{\theta(4H-1)}{(2H)^2}\bigg\{1+\frac{\Gamma(3-4H)\Gamma(4H-1)}{\Gamma(2H)\Gamma(2-2H)}\bigg\}.
\eea
\halflineskip

\begin{lemma}\label{202307061023}
Let $m\geq2$. Assume $H\in \left( \frac{1}{2}, \frac{m+1}{2m}\right)$. Then 
$C_U(m,H,\theta)$ is finite and 
\bea\label{202307061021}
E\big[\Gamma^{(m)}(U_T,...,U_T)\big]
&=&
2^{m-1}\big\langle \underbrace{u_T\otimes_1\cdots\otimes_1u_T}_{m-1},u_T\big\rangle_{\calh^{\otimes2}}
\nn\\&=& 
T^{-\half(m-2)}{\vred C_U(m,H,\theta)}
+o\big(T^{-\half(m-2)}\big)
\eea
as $T\to\infty$. 
\end{lemma}
\proof 
Let
\bea\label{202307080224}
I_{T}
&=&
\int_{[0, T] ^{2m}}a(x_{1}, x_{2}, x_{3})a(x_{3}, x_{4}, x_{5})\ldots a(x_{2m-1}, x_{2m}, x_{1})dx_{1}\cdots dx_{2m}
\eea
and 
\beas
I'_{\infty}
&=&
2m\int_{(0,\infty)^{2m-1}}
a(0, x_{2}, x_{3}) a(x_{3}, x_{4}, x_{5})\cdots a(x_{2m-3},x_{2m-2}, x_{2m-1})a(x_{2m-1}, x_{2m},0)\>dx_{2}\cdots dx_{2m}. 
\eeas

From (\ref{202307060939}), we obtain 
\bea\label{202308091838}
(2m)^{-1}I'_\infty
&\leq&
\int_{\bbR^{m-1}}\ol{a}(x_1)\ol{a}(x_1-x_2)\cdots
\ol{a}(x_{m-2}-x_{m-1})\ol{a}(x_{m-1})dx_1\cdots dx_{m-1}, 
\eea
and $I'_\infty<\infty$ by using the estimate (\ref{202307050606}) of Lemma \ref{202307060258}, and 
Lemma \ref{202308120413}.

By L'H\^opital's rule, 
\bea\label{202307061019}
 \lim_{T\to \infty} \frac{ I_{T}}{T} 
 &=& 
  \lim_{T\to \infty} \frac{dI_T}{dT}
 \nn\\&=&
 2m \lim_{T\to \infty}\int_{[0, T] ^{2m-1}} 
 a (T, x_{2}, x_{3})a (x_{3}, x_{4}, x_{5})\cdots a(x_{2m-1}, x_{2m}, T)dx_{2}...dx_{2m}
 \nn\\&=& 
 2m \lim_{T\to \infty}\int_{[0, T] ^{2m-1}}
a(0, x_{2}, x_{3})a (x_{3}, x_{4}, x_{5})\cdots a(x_{2m-1}, x_{2m}, 0) dx_{2}...dx_{2m}
\nn\\&=&
I'_\infty,
\eea
where we changed variables as $\tilde{x}_{i}= T-x_{i}$ for $i=2,...,m$. 

From (\ref{202307030144}) {\sblue and the expression of the scalar product in $\calh^{\otimes2}$}, 
\bea\label{202307061020}
E\big[\Gamma^{(m)}(U_T,...,U_T)\big]
&= &
2^{m-1}\big\langle \underbrace{u_T\otimes_1\cdots\otimes_1u_T}_{m-1},u_T\big\rangle_{\calh^{\otimes2}}
\nn\\&=&
{\vred 2^{m-1}K_U^mT^{-m/2}\alpha_H^m}I_T
\eea
for $m\geq2$. 
Now we obtain (\ref{202307061021}) from (\ref{202307061019}) and (\ref{202307061020})
since ${\vred C_U(m,H,\theta)
=2^{m-1}K_U^m\alpha_H^mI'_\infty}$.
\qed\halflineskip
\begin{lemma}\label{202307080000}
Let $m\geq2$. Suppose that $H=\frac{m+1}{2m}$. Then, for any $\ep>0$, 
\bea\label{202307080001}
E\big[\Gamma^{(m)}(U_T,...,U_T)\big]
\yeq
2^{m-1}\big\langle \underbrace{u_T\otimes_1\cdots\otimes_1u_T}_{m-1},u_T\big\rangle_{\calh^{\otimes2}}
\yeq
{\sblue o}\big(T^{-\half(m-2)+\ep}\big)
\eea
as $T\to\infty$. 
\end{lemma}
\proof
Recall that the functions $a_T(x,z,y)$, $\ol{a}(x)$ are associated with $H=\frac{m+1}{2m}$. 
%
By (\ref{202307080224}) and (\ref{202307060939}), 
\bea\label{202307080228}
I_{T}
&\leq&
\int_{[0, T]^m}\ol{a}(x_1-x_2)\ol{a}(x_2-x_3)\cdots\ol{a}(x_{m-1}-x_1)dx_1\cdots dx_m.
\eea

For any $\ep_1>0$, 
Lemma \ref{202307060258} yields 
\bea\label{202307080251}
\ol{a}(r)
&\leq& 
C(1\wedge r^{2H-2})
\yeq
C\big(1\wedge r^{2H-2-\ep_1}(r/T)^{\ep_1}T^{\ep_1}\big)
\nn\\&\leq& 
\wt{a}(r)T^{\ep_1}
\quad (\forall r\in(0,T);\>T\geq1), 
\eea
where 
$\wt{a}(x)=C\big(1\wedge |x|^{2H-2-\ep_1}\big)$ for $x\in\bbR$. 
Let 
\bea\label{202307080249}
\wt{I}_{T}
&=&
\int_{[0, T]^m}\wt{a}(x_1-x_2)\wt{a}(x_2-x_3)\cdots\wt{a}(x_{m-2}-x_{m-1})\wt{a}(x_{m-1}-x_1)dx_1\cdots dx_m. 
\eea
Then
\bea\label{202308091846}
\lim_{T\to\infty}\frac{d\wt{I}_T}{dT}
&=&
m\lim_{T\to\infty}\int_{[0, T]^{m-1}}\wt{a}(T-x_2)\wt{a}(x_2-x_3)\cdots\wt{a}(x_{m-2}-x_{m-1})\wt{a}(x_{m-1}-T)dx_2\cdots dx_m
\nn\\&=&
m\lim_{T\to\infty}\int_{[0, T]^{m-1}}\wt{a}(x_2)\wt{a}(x_2-x_3)\cdots\wt{a}(x_{m-2}-x_{m-1})\wt{a}(x_{m-1})dx_2\cdots dx_m
\quad(x_i\leftarrow T-x_i)
\nn\\&=&
m\int_{[0,\infty)^{m-1}}\wt{a}(x_2)\wt{a}(x_2-x_3)\cdots\wt{a}(x_{m-2}-x_{m-1})\wt{a}(x_{m-1})dx_2\cdots dx_m
\>=:\>\wt{I}_\infty'
\eea
The limit 
$\wt{I}_\infty'$ is finite by Lemma \ref{202308120413} applied to $\alpha_i(x)=\wt{a}(x)=C\big(1\wedge |x|^{2H-2-\ep_1}\big)$. 
\begin{en-text}
since 
$\wt{I}_\infty'\leq
m\big\|\wt{a}^{*(m-1)}\wt{a}\big\|_{L^1(\bbR)}
\yleq 
m\|\wt{a}\|_{L^{\frac{m}{m-1}}(\bbR)}^m<\infty$ and 
$(2H-2-\ep_1)\frac{m}{m-1}=-1-\frac{\ep_1m}{m-1}<-1$. 
\end{en-text}

Set $\ep_1=\ep/m$ for a given $\ep>0$. 
Now, (\ref{202307080228}) and (\ref{202307080251}) give 
$
I_{T}
\leq
T^{m\ep_1}\wt{I}_{T}$.
Therefore, from (\ref{202307061020}), 
\beas
0
&\leq&
T^{\half(m-2)-\ep}E\big[\Gamma^{(m)}(U_T,...,U_T)\big]
\yeq
{\vred 2^{m-1}K_U^mT^{-m/2}\alpha_H^m}\> T^{-1-\ep}I_T
\nn\\&\leq&
{\vred 2^{m-1}K_U^mT^{-m/2}\alpha_H^m}\> T^{-1}\wt{I}_T
\>\to_{T\to\infty}\>
{\vred 2^{m-1}K_U^mT^{-m/2}\alpha_H^m}\wt{I}_\infty'\><\>\infty
\eeas
by L'H\^opital's rule. This completes the proof. 
\qed\halflineskip

For $p_1,...,p_m\in\bbR$, define $\sfB_m(p_1,p_2,...,p_m)$ 
by
\beas 
\sfB_m(p_1,p_2,...,p_m)
&=&
\int_{[0,1]^m}|x_1-x_2|^{p_1}|x_2-x_3|^{p_2}\cdots|x_{m-1}-x_m|^{p_{m-1}}|x_m-x_1|^{p_m}dx_1...dx_m\in[0,\infty].
\eeas
Define $a_T(x,y)$ by 
\bea\label{202308091907}
a_T(x,y) 
&=& 
\int_0^Ta(x,z,y)dz
\yeq 
\int_0^T e^{-\theta|x-z|}|z-y|^{2H-2}dz
\quad(x,y\in\bbR). 
\eea
\begin{lemma}\label{202308120342}
Let $m\geq2$. 
Suppose that $H\in(\frac{m+1}{2m},1)$. Then ${\sf B}_m{\bred(2H-2,...,2H-2)}<\infty$ and 
\beas&&
\lim_{T\to\infty}
T^{-(2H-1)m}\int_{[0,T]^{m}}a_T(x_1,x_2)a_T(x_2,x_3)\ldots a_T(x_m,x_1)dx_{1}...dx_{m}
\nn\\&&\hspace{30pt}\yeq
2^m\theta^{-m}{\sf B}_m{\sblue(2H-2,...,2H-2)}. 
\eeas
\end{lemma}
\proof
We have 
\bea\label{202306220818}
a_T(x,y)
&=& 
2T^{2H-2}A_T(T^{-1}x,T^{-1}y),
\eea
where 
\beas 
A_T(x,y) 
&=& 
\frac{T}{2}\int_0^1e^{-T\theta|x-z|}|z-y|^{2H-2}dz.
\eeas
By (\ref{202307060939}) and Lemma \ref{202307060258}, for some constant $C$, 
$a_T(x,y)\leq C|x-y|^{2H-2}$ for $x,y\in\bbR$, 
in particular, 
\bea\label{202306230159}
A_T(x,y)
\yeq
2^{-1}T^{-2H+2}a_T(Tx,Ty)
\yleq
2^{-1}T^{-2H+2}\ol{a}(|Tx-Ty|)
\yleq
2^{-1}C|x-y|^{2H-2}
\eea
for $x,y\in\bbR$. 
\begin{en-text}
\beas 
0\leq A_T(x,y) 
&=&
\frac{T}{2}\int_{[0,1]\cap\{|z-y|\geq2^{-1}|x-y|\}}e^{-T\theta|x-z|}|z-y|^{2H-2}dz
\nn\\&&
+\frac{T}{2}\int_{[0,1]\cap\{|z-y|<2^{-1}|x-y|\}}e^{-T\theta|x-z|}|z-y|^{2H-2}dz
\nn\\&\leq& 
2^{1-2H}|x-y|^{2H-2}T\int_{[0,1]\cap\{|z-y|\geq2^{-1}|x-y|\}}e^{-T\theta|x-z|}dz
\nn\\&&
+Te^{-2^{-1}T\theta|x-y|}\int_{[0,1]\cap\{|z-y|<2^{-1}|x-y|\}}|z-y|^{2H-2}dz
\eeas
\end{en-text}
Furthermore, 
{\sblue by using the convergence of the Laplace distribution to the delta-measure,} 
it is not difficult to show 
\bea\label{202306230203}
A_T(x,y) &\to& \theta^{-1}|x-y|^{2H-2}\quad(T\to\infty)
\eea
for $(x,y)\in(0,1)^2$, $x\not=y$. 
Lebesgue's theorem with (\ref{202306230159}) and (\ref{202306230203}) ensures 
\beas&&
T^{-(2H-1)m}\int_{[0,T]^{m}}a_T(x_1,x_2)a_T(x_2,x_3)\ldots a_T(x_m,x_1)dx_{1}...dx_{m}
\nn\\&=&
2^m\int_{[0,1]^{m}}A_T(x_1,x_2)A_T(x_2,x_3)\ldots A_T(x_m,x_1)dx_{1}...dx_{m}
\nn\\&\to&
2^m\theta^{-m}{\sf B}_m{\sblue(2H-2,...,2H-2)}\quad(T\to\infty)
\eeas
if $\sfB_m(2H-2,...,2H-2)<\infty$. 
However, we know $\sfB_m(2H-2,...,2H-2)<\infty$ when $H>\frac{m+1}{2m}$. 
{\sblue See Lemma \ref{202402110703} below.} 
\qed\halflineskip

{\sblue
%
\begin{en-text}
Recall the definition of the function 
$\sfB_m(p_1,p_2,...,p_m)$ taking values in $[0,\infty]$: 
\beas 
\sfB_m(p_1,p_2,...,p_m)
&=&
\int_{[0,1]^m}|x_1-x_2|^{p_1}|x_2-x_3|^{p_2}\cdots|x_{m-1}-x_m|^{p_{m-1}}|x_m-x_1|^{p_m}dx_1...dx_m\in[0,\infty].
\eeas
for $p_1,...,p_m\in\bbR$, $m\in\bbZ_{\geq2}$. 
\end{en-text}
\begin{lemma}\label{202402110703}
Let $m\in\bbZ_{\geq2}$. 
Suppose that the numbers $p_1,...,p_m>-1$ satisfy $\sum_{i=1}^mp_i+m-1>0$. Then 
$\sfB_m(p_1,p_2,...,p_m)<\infty$. 
\end{lemma}
\proof
The variance gamma distribution $\text{VG}(\lambda,\alpha,\beta,\mu)$ is a probability distribution 
on $\bbR$ with the density function
\beas
p(x) 
&=& 
\frac{1}{\sqrt{\pi}\Gamma(\lambda)}(\alpha^2-\beta^2)^\lambda
\left(\frac{|x-\mu|}{2\alpha}\right)^{\lambda-\half}K_{\lambda-\half}(\alpha|x-\mu|)\exp(\beta(x-\mu))
\quad(x\in\bbR), 
\eeas
where $\lambda,\alpha\in(0,\infty)$, $\beta\in\bbR$ ($\alpha>|\beta|$) and $\mu\in\bbR$ are parameters, 
and $K_\nu$ is the Bessel function of the third kind with index $\nu$.  
See e.g. Iacus and Yoshida \cite{iacus2017simulation} for the variance gamma distribution 
and the related variance gamma process. 
Here we will use the variance gamma distribution 
$\text{VG}(\lambda,1,0,0)$ for $\lambda>0$. 
Denote the density 
of $\text{VG}(\lambda,1,0,0)$ by $p(x;\lambda)$. 

The following facts are known: 
\begin{enumerate}[(i)] 
\im$K_{-\nu}(z)=K_\nu(z)$
\im $K_\nu(z)\sim 2^{-1}\Gamma(\nu)(z/2)^{-\nu}$ as $z\to0$ when $\text{Re}(\nu)>0$, 
and $K_0(z)\sim-\log z$. 
\im As $z\to\infty$ under $|\arg z|\leq3\pi/2-\ep$ with $\ep>0$, 
\beas 
K_\nu(z) 
&\sim& 
\sqrt{\frac{\pi}{2z}}e^{-z}\sum_{k=0}^\infty\frac{a_k(\nu)}{z^k}, 
\eeas
where 
\beas 
a_k(\nu) 
&=& 
\frac{(4\nu^2-1)(4\nu^2-3^2)\cdots(4\nu^2-(2k-1)^2)}{8^kk!}.
\eeas
\end{enumerate}
Around $x=0$, 
the density function $p(x;\lambda)$ has the singularity $|x|^{2\lambda-1}$ when $2\lambda-1<0$, 
$-\log |x|$ when $2\lambda-1=0$, and no singularity when $2\lambda-1>0$. 
Moreover, the function $p(x;\lambda)$ rapidly decays when $|x|\to\infty$. 
Thus, we have the estimate 
\bea\label{202402110829}
|x|^{2\lambda-1}1_{\{|x|\leq1\}}&\leq& C_\lambda\> p(x;\lambda)\qquad(x\in\bbR)
\eea
for some constant $C_\lambda$ depending on $\lambda>0$. 

The family of variance gamma distributions is closed under convolution. 
In fact, in our case, the characteristic function of $\text{VG}(\lambda,1,0,0)$ is 
\beas
\varphi_{\text{VG}(\lambda,1,0,0)}(u) 
&=& 
\big(1+u^2\big)^{-\lambda}
\quad(u\in\bbR)
\eeas
and hence 
\bea\label{202402110830}
\text{VG}(\lambda_1,1,0,0)*\text{VG}(\lambda_2,1,0,0)
&=& 
\text{VG}(\lambda_1+\lambda_2,1,0,0)
\eea
for $\lambda_1,\lambda_1>0$. 

Suppose that $p_i>-1$ for $i=1,...,m$. 
Let $\lambda_i=(p_i+1)/2>0$ for $i=1,...,m$. 
Then 
\beas
&&
\int_{[0,1]^m}|x_1-x_2|^{p_1}|x_2-x_3|^{p_m}\cdots
|x_{m-1}-x_m|^{p_{m-1}}|x_m-x_1|^{p_m}
dx_1...dx_m
\nn\\&\simleq&
\int_{\bbR^m} 1_{[0,1]}(x_1)p(x_1-x_2;\lambda_1)p(x_2-x_3;\lambda_2)\cdots 
p(x_{m-1}-x_m;\lambda_{m-1})p(x_m-x_1;\lambda_m)dx_1...dx_m
\quad((\ref{202402110829}))
\nn\\&=&
\int_{\bbR^2} 1_{[0,1]}(x_1)p(x_1-x_m;\lambda_1+\cdots+\lambda_{m-1})p(x_m-x_1;\lambda_m)dx_mdx_1
\quad((\ref{202402110830}))
\nn\\&=&
\int_{\bbR} 1_{[0,1]}(x_1)p(0;\lambda_1+\cdots+\lambda_{m})dx_1
\quad((\ref{202402110830}))
\nn\\&=&
p(0;\lambda_1+\cdots+\lambda_{m}). 
\eeas
On the other hand, 
$p(0;\lambda_1+\cdots+\lambda_{m})<\infty$ since 
the density function $p(x;\lambda_1+\cdots+\lambda_{m})$ has no singularity at the origin 
due to 
\beas 
2(\lambda_1+\cdots+\lambda_{m})-1
\yeq
\sum_{i=1}^mp_i+m-1
\>>\>0
\eeas
by assumption. 
\qed\halflineskip

Under the assumption of Lemma \ref{202308120342}, obviously 
Lemma \ref{202402110703} ensures $\sfB_m(2H-2,...,2H-2)<\infty$ 
since $2H-2>-1$ by $H>1/2$, and 
$m(2H-2)+m-1=2mH-m-1>0$. 
}
\halflineskip

\begin{lemma}\label{202307061046}
Let $m\geq2$ and ${\vred C_U'(m,H,\theta)=2^{2m-1}K_U^m\alpha_H^m\theta^{-m}}\sfB_m(2H-2,...,2H-2)$. 
Suppose that $H\in(\frac{m+1}{2m},1)$. Then $C_U'(m,H,\theta)<\infty$ and 
\bea\label{202306230217}
T^{(\frac{3}{2}-2H)m}
E\big[\Gamma^{(m)}(U_T,...,U_T)\big]
&=&
2^{m-1}
\big\langle \underbrace{u_T\otimes_1\cdots\otimes_1u_T}_{m-1},u_T\big\rangle_{\calh^{\otimes2}}T^{(\frac{3}{2}-2H)m}
\nn\\&\to&
{\vred C_U'(m,H,\theta)}
\eea
as $T\to\infty$. 
\end{lemma}
\proof 
From (\ref{202307030127}) and (\ref{202306220818}), we obtain 
\bea\label{202307080045}&&
E\big[\Gamma^{(m)}(U_T,...,U_T)\big]
\nn\\&=& 
2^{m-1}\big\langle u_T\otimes_1\cdots\otimes_1u_T,u_T\big\rangle_{\calh^{\otimes2}}
\nn\\&=&
{\vred2^{m-1}K_U^m\alpha_H^m}T^{-m/2}\int_{[0,T]^{m}}a_T(x_1,x_2)a_T(x_2,x_3)\ldots a_T(x_m,x_1)dx_{1}...dx_{m}
. 
\eea
\begin{en-text}
Since 
\beas&&
T^{-m(2H-2)-2^{-1}m}E\big[\Gamma^{(m)}(U_T,...,U_T)\big]
\nn\\&=&
C_U'(m,H,\theta,\sigma)\theta^m\int_{[0,1]^{m}}A_T(x_1,x_2)A_T(x_2,x_3)\ldots A_T(x_m,x_1)dx_{1}...dx_{m}, 
\eeas
by (\ref{202306220818}), 
\end{en-text}
Now the convergence (\ref{202306230217}) follows from Lemma \ref{202308120342}. 
\qed\halflineskip

\subsection{Expansion of $E\big[\Gamma^{(2)}(U_T,U_T)\big]$}
Let 
\bea\label{202308150222}
{\vred C_U''(2,H,\theta)}
&=&
{\bred-}
{\vred 
\frac{(2H-1)\theta^{4H-2}}{2H^2(3-4H)\Gamma(2H)^2}.
}
\eea
\begin{lemma}\label{202307081346}
Suppose that $H\in(1/2,3/4)$. Then 
\beas
E\big[\Gamma^{(2)}(U_T,U_T)\big]
&=&
2\big\langle u_T,u_T\big\rangle_{\calh^{\otimes2}}
\nn\\&=&
C_U(2,H,\theta)
+{\vred C_U''(2,H,\theta)}T^{4H-3}+o(T^{4H-3})
\eeas
as $T\to\infty$.
\end{lemma}
\proof
From (\ref{202307061020}), 
\bea\label{202307081341}
E\big[\Gamma^{(2)}(U_T,U_T)\big]
&= &
2\big\langle u_T,u_T\big\rangle_{\calh^{\otimes2}}
\yeq
{\vred2K_U^2}\alpha_H^2T^{-1}I^{(2)}_T,
\eea
where 
\beas
I^{(2)}_{T}
&=&
\int_{[0, T] ^{4}}a(x_{1}, x_{2}, x_{3})a(x_{3}, x_{4}, x_{1})dx_{1}\cdots dx_{4}.
\eeas

In Lemma \ref{202307061023} and its proof, we already know 
\begin{en-text}
\beas
I^{(2)\prime}_\infty
&:=& 
\lim_{T\to\infty}\frac{I^{(2)}_T}{T}
\yeq 
\lim_{T\to\infty}\frac{dI^{(2)}_T}{dT}
\nn\\&=&
\lim_{T\to\infty}4\int_{[0,T]^3}a(0,x_2,x_3)a(x_3,x_4,0)dx_2dx_3dx_4
\nn\\&=&
4\int_{[0,\infty)^3}a(0,x_2,x_3)a(x_3,x_4,0)dx_2dx_3dx_4. 
\eeas
\end{en-text}
\beas
\frac{dI^{(2)}_T}{dT}
&=&
4\int_{[0,T]^3}a(0,x_2,x_3)a(x_3,x_4,0)dx_2dx_3dx_4
\eeas
and
\bea\label{202307081342}
I^{(2)\prime}_\infty
&:=& 
\lim_{T\to\infty}\frac{I^{(2)}_T}{T}
\yeq 
\lim_{T\to\infty}\frac{dI^{(2)}_T}{dT}
\yeq
\big({\vred 2K_U^2\alpha_H^2}\big)^{-1}C_U(2,H,\theta). 
\eea

\begin{en-text}
Since, as will be shown, the limit in the rightmost expression of (\ref{202307081320}) below exists, by L'H\^opital's rule, we obtain 
\bea\label{202307081320}
\lim_{T\to\infty}\big(T^{-1}I^{(2)}_T-I^{(2)\prime}_\infty\big)/T^{4H-3}
&=&
\lim_{T\to\infty}\bigg(T^{-1}\frac{dI^{(2)}_T}{dT}-T^{-2}I^{(2)}_T\bigg)/\big((4H-3)T^{4H-4}\big)
\nn\\&=&
\lim_{T\to\infty}\bigg(\frac{dI^{(2)}_T}{dT}-T^{-1}I^{(2)}_T\bigg)/\big((4H-3)T^{4H-3}\big)
\nn\\&=&
\lim_{T\to\infty}\bigg(\frac{d^2I^{(2)}_T}{dT^2}-T^{-1}\frac{dI^{(2)}_T}{dT}+T^{-2}I^{(2)}_T\bigg)/\big((4H-3)^2T^{4H-4}\big)
\nn\\&=&
\lim_{T\to\infty}\frac{dI^{(2)}_T}{dT}/\big((4H-3)T^{4H-1}\big)
\nn\\&&\text{(The equality is valid if this limit exists.}
\nn\\&&\text{L'H\^opital's rule can apply since the numerator diverges.)}
\nn\\&=&
\lim_{T\to\infty}\big(\frac{dI^{(2)}_T}{dT}-I^{(2)\prime}_\infty\big)/\big((4H-2)T^{4H-3}\big)
\nn\\&=&
\lim_{T\to\infty}\frac{dI^{(2)}_T}{dT}\big(I^{(2)}_T-I^{(2)\prime}_\infty\big)/\big((4H-3)T^{4H-4}\big)
\nn\\&=&
\lim_{T\to\infty}4(3-4H)^{-1}T^{4-4H}\big(I^{(2,1)}_T+I^{(2,2)}_T+I^{(2,3)}_T\big), 
\eea
\end{en-text}

In the following equalities of (\ref{202307081320}), 
$=^{***}$ is obvious, and $=^{**}$ is verified by L'H\^opital's rule with the aid of $\frac{dI^{(2)}_T}{dT}-I^{(2)\prime}_\infty\to0$ {\sblue as $T\to\infty$}. 
As will be seen, the limit on the right-hand side of $=^{***}$ is {\bred non-zero}. 
Therefore, $I^{(2)}_T-TI^{(2)\prime}_\infty=\int_0^T \big(\frac{dI^{(2)}_t}{dt}-I^{(2)\prime}_\infty\big)dt
\to\infty$ since $\int_1^\infty t^{4H-3}dt=\infty$. 
With this fact, L'H\^opital's rule applies to the equalities $=^*$.
In this way, we obtain
\bea\label{202307081320}
\lim_{T\to\infty}\big(T^{-1}I^{(2)}_T-I^{(2)\prime}_\infty\big)/T^{4H-3}
&=&
\lim_{T\to\infty}\big(I^{(2)}_T-TI^{(2)\prime}_\infty\big)/T^{4H-2}
\nn\\&=^*&
\lim_{T\to\infty}\big(\frac{dI^{(2)}_T}{dT}-I^{(2)\prime}_\infty\big)/\big((4H-2)T^{4H-3}\big)
\nn\\&=^{**}&
\lim_{T\to\infty}\frac{d^2I^{(2)}_T}{dT^2}/\big((4H-2)(4H-3)T^{4H-4}\big)
\nn\\&=^{***}&
\lim_{T\to\infty}4(4H-2)^{-1}({\bred4H-3})^{-1}T^{4-4H}\big(I^{(2,1)}_T+I^{(2,2)}_T+I^{(2,3)}_T\big), 
\nn\\&&
\eea
where 
\beas
I^{(2,1)}_T
&=&
\int_{[0,T]^2}a(0,T,x_3)a(x_3,x_4,0)dx_3dx_4, 
\eeas
\beas
I^{(2,2)}_T
&=&
\int_{[0,T]^2}a(0,x_2,T)a(T,x_4,0)dx_2dx_4
\eeas
and 
\beas
I^{(2,3)}_T
&=&
\int_{[0,T]^2}a(0,x_2,x_3)a(x_3,T,0)dx_2dx_3. 
\eeas

For $I^{(2,i)}_T$ ($i=1,2,3$), we have the following estimates: 
\bea\label{202307081321}
I^{(2,1)}_T
&=&
\int_{[0,T]^2}e^{-\theta T}|T-x_3|^{2H-2}e^{-\theta|x_3-x_4|}|x_4|^{2H-2}dx_3dx_4
\>\simleq\>
e^{-\theta T/2},
\eea
\bea\label{202307081322}
I^{(2,2)}_T
&=&
\int_{[0,T]^2}e^{-\theta|x_2|}|x_2-T|^{2H-2}e^{-\theta|T-x_4|}|x_4|^{2H-2}dx_2dx_4
\nn\\&=&
T^2\int_{[0,1]^2}e^{-\theta Tx_2}|Tx_2-T|^{2H-2}e^{-\theta|T-Tx_4|}|Tx_4|^{2H-2}dx_2dx_4
\nn\\&=&
T^{4H-2}\int_{[0,1]^2}e^{-\theta Tx_2}|x_2-1|^{2H-2}e^{-\theta T|1-x_4|}|x_4|^{2H-2}dx_2dx_4
\nn\\&=&
T^{4H-4}{\colorr\theta^{-2}}\int_{[0,1]^2}\theta Te^{-\theta Tx_2}|x_2-1|^{2H-2}\ \theta Te^{-\theta T|1-x_4|}|x_4|^{2H-2}dx_2dx_4
\nn\\&\sim&
T^{4H-4}{\colorr\theta^{-2}}
\eea
and 
\bea\label{202307081323}
I^{(2,3)}_T
&=&
\int_{[0,T]^2}e^{-\theta x_2}|x_2-x_3|^{2H-2}e^{-\theta|x_3-T|}T^{2H-2}dx_2dx_3
\nn\\&=&
T^{2H}\int_{[0,1]^2}e^{-\theta Tx_2}|Tx_2-Tx_3|^{2H-2}e^{-\theta|Tx_3-T|}dx_2dx_3
\nn\\&=&
T^{4H-2}\int_{[0,1]^2}e^{-\theta Tx_2}|x_2-x_3|^{2H-2}e^{-\theta T|1-x_3|}dx_2dx_3
\nn\\&=&
T^{4H-4}{\colorr\theta^{-2}}\int_{[0,1]^2}\theta Te^{-\theta Tx_2}|x_2-x_3|^{2H-2}\ \theta Te^{-\theta T|1-x_3|}dx_2dx_3
\nn\\&\sim&
T^{4H-4}{\colorr\theta^{-2}}
\eea
as $T\to\infty$. 

Thus, we obtain 
\bea\label{202307081343}
\lim_{T\to\infty}\big(T^{-1}I^{(2)}_T-I^{(2)\prime}_\infty\big)/T^{4H-3}
&=&
\lim_{T\to\infty}4{\vred(4H-2)^{-1}}({\bred4H-3})^{-1}T^{4-4H}\big(I^{(2,1)}_T+I^{(2,2)}_T+I^{(2,3)}_T\big)
\nn\\&=&
{\colorr8}(4H-2)^{-1}({\bred4H-3})^{-1}\theta^{-2}
\eea
as $T\to\infty$, 
from (\ref{202307081320}), (\ref{202307081321}), (\ref{202307081322}) and (\ref{202307081323}). 

From (\ref{202307081341}),  (\ref{202307081342}) and (\ref{202307081343}), 
\beas
E\big[\Gamma^{(2)}(U_T,U_T)\big]
&=&
{\vred2K_U^2}\alpha_H^2T^{-1}I^{(2)}_T
\nn\\&=&
C_U(2,H,\theta)
+
{\vred C_U''(2,H,\theta)}
T^{4H-3}
+o(T^{4H-3})
\eeas
as $T\to\infty$. 
This completes the proof. 
\qed\halflineskip

\subsection{Estimate of $U_T$, $V_T$ and $W_T$}
{\sblue The $(s,p)$-Sobolev norm of functional $F$ is defined as $\|F\|_{s,p}=\|(1-L)^{s/2}F\|_p$ 
for $s\in\bbR$ and $p>1$. }
{\ared Let $D_\infty=\cap_{s\in\bbR,p>1}D_{s,p}$.}
\begin{lemma}\label{20230813021}
$U_T=O_{D_\infty}(1)$, i.e., 
$\|U_T\|_{s,p}=O(1)$ as $T\to\infty$ for every $s\in\bbR$ and $p>1$. 
\end{lemma}
\proof $
E[U_T^2]=2\langle u_T,u_T\rangle_{\mfh^{\otimes2}}=E\big[\Gamma^{(2)}(U_T,U_T)\big]
=O(1)
$ thanks to Lemma \ref{202307081346}. 
Hypercontractivity and a fix chaos give the result. 
\qed\halflineskip
\begin{lemma}\label{202308110225}
$V_T=O_{D_\infty}(T^{-1/2})$. 
\end{lemma}
\proof 
We have 
\beas
E[V_T^2] 
&=& 
2\langle v_T,v_T\rangle_{\mfh^{\otimes2}}
\nn\\&=&
2\alpha_H^2K_V^2T^{-1}\int_{[0,T]^4} e^{-\theta (T-t_1)-\theta (T-t_2)}|t_2-t_3|^{2H-2}e^{-\theta (T-t_3)-\theta (T-t_4)}
|t_4-t_1|^{2H-2}dt_1dt_2dt_3dt_4
\nn\\&\simleq&
T^{-1}\int_{[0,T]^2} e^{-\theta (T-t_1)}|T-t_3|^{2H-2}e^{-\theta (T-t_3)}
|T-t_1|^{2H-2}dt_1dt_3
\nn\\&&
\quad{\sblue(\text{Use (\ref{202307060939}) and (\ref{202307050606}) for 
the integrals with respect to $t_2$ and $t_4$})}
\nn\\&\leq&
T^{-1}\bigg(\int_{[0,\infty)}e^{-\theta t}t^{2H-2}dt\bigg)^2
\yeq \big(T^{{\sblue -1/2}}\theta^{1-2H}\Gamma(2H-1)\big)^2
\eeas
for all $T>0$. 
Then we obtain the results by hypercontractivity. 
\qed\halflineskip
\begin{lemma}\label{20230813034}
$W_T=O_{D_\infty}(T^{-1/2})$. 
\end{lemma}
\proof 
It is sufficient to observe that 
\beas 
E[W_T^2] 
&=& 
\langle w_T,w_T\rangle_\mfh
\nn\\&=&
T^{-1}K_W^2\alpha_H\int_{[0,T]^2}(e^{-\theta t}-e^{-2\theta T+\theta t})|t-s|^{2H-2}
(e^{-\theta s}-e^{-2\theta T+\theta s})dtds
\nn\\&\leq&
T^{-1}K_W^2\alpha_H\int_{[0,T]^2}e^{-\theta t}|t-s|^{2H-2}e^{-\theta s}dtds
\nn\\&&
\quad(\because 0\leq e^{-\theta t}-e^{-2\theta T+\theta t}=e^{-\theta t}(1-e^{-2\theta (T-t)})\leq e^{-\theta t})
\nn\\&\simleq&
T^{-1}\int_{[0,T]}e^{-\theta t}t^{2H-2}dt\quad(\text{(\ref{202307060939}) and (\ref{202307050606}))}
\nn\\&\leq&
T^{-1}\theta^{1-2H}\Gamma(2H-1)
\eeas
for all $T>0$. 
\qed

\subsection{Cross-gamma factors}

\begin{en-text}
\begin{lemma}\label{202308091510}
$
E\big[\Gamma^{(2)}(U_T,V_T)\big]
=
E\big[\Gamma^{(2)}(V_T,U_T)\big]
=
O(T^{-1})
$ as $T\to\infty$. 
\end{lemma}
\proof 
We have 
\beas 
E\big[\Gamma^{(2)}(U_T,V_T)\big]
&=& 
E\big[\Gamma^{(2)}(V_T,U_T)\big]
\yeq 
2^{-1}E\big[\langle DU_T,DV_T\rangle\big]
\nn\\&=&
2\langle u_T,v_T\rangle_{\mfh^{\otimes2}}
\yeq
C(\theta,H)T^{-1}J_T,
\eeas
where $C(\theta,H)$ is a constant and 
\beas 
J_T 
&=& 
\int_{[0,T]^4}e^{-\theta|t_1-s_1|}|s_1-s_2|^{2H-2}e^{-\theta|T-s_2|-\theta|T-t_2|}|t_2-t_1|^{2H-2}ds_1ds_2dt_1dt_2. 
\eeas
Then we have 
\bea\label{202308091449}
J_T &=& O(1)
\eea
as $T\to\infty$. 
Indeed, by using (\ref{202307060939}), and (\ref{202307050606}) of Lemma \ref{202307060258}, 
we obtain 
\beas 
J_T 
&\simleq&
\int_{[0,T]^2}\big(1\wedge|t_1-s_2|^{2H-2}\big)e^{-\theta(T-s_2)}\big(1\wedge|T-t_1|^{2H-2}\big)ds_2dt_1
\nn\\&\simleq&
\int_{[0,T]}\big(1\wedge|T-t_1|^{2H-2}\big)\big(1\wedge|T-t_1|^{2H-2}\big)dt_1
\nn\\&\leq&
\int_{[0,T]}\big(1\wedge|T-t_1|^{4H-4}\big)dt_1
\yeq 
\int_{[0,\infty)}\big(1\wedge t^{4H-4}\big)dt
\><\>\infty
\eeas
due to $4H-4<-1$ when $H<3/4$. 
This completes the proof. 
\qed\halflineskip
\end{en-text}
\begin{lemma}\label{202308091510a}
$
E\big[\Gamma^{(2)}(U_T,V_T)\big]
=
E\big[\Gamma^{(2)}(V_T,U_T)\big]
=
O(T^{-1})
$ as $T\to\infty$. 
\end{lemma}
\proof 
We have 
\beas 
E\big[\Gamma^{(2)}(U_T,V_T)\big]
&=& 
E\big[\Gamma^{(2)}(V_T,U_T)\big]
\yeq 
2^{-1}E\big[\langle DU_T,DV_T\rangle\big]
\nn\\&=&
2\langle u_T,v_T\rangle_{\mfh^{\otimes2}}
\yeq
C(\theta,H)T^{-1}J_T,
\eeas
where $C(\theta,H)$ is a constant and 
\beas 
J_T 
&=& 
\int_{[0,T]^4}e^{-\theta|t_1-s_1|}|s_1-s_2|^{2H-2}e^{-\theta|T-s_2|-\theta|T-t_2|}|t_2-t_1|^{2H-2}ds_1ds_2dt_1dt_2. 
\eeas
Then we have 
\bea\label{202308091449}
J_T &=& O(1)
\eea
as $T\to\infty$. 
Indeed, by using (\ref{202307060939}), and (\ref{202307050606}) of Lemma \ref{202307060258}, 
we obtain 
\beas 
J_T 
&\simleq&
\int_{[0,T]^2}\big(1\wedge|t_1-s_2|^{2H-2}\big)e^{-\theta(T-s_2)}\big(1\wedge|T-t_1|^{2H-2}\big)ds_2dt_1
\nn\\&\simleq&
\int_{[0,T]}\big(1\wedge|T-t_1|^{2H-2}\big)\big(1\wedge|T-t_1|^{2H-2}\big)dt_1
\nn\\&\leq&
\int_{[0,T]}\big(1\wedge|T-t_1|^{4H-4}\big)dt_1
{\sblue\><\>}
\int_{[0,\infty)}\big(1\wedge t^{4H-4}\big)dt
\ <\ \infty
\eeas
due to $4H-4<-1$ when $H<3/4$. 
\qed\halflineskip

\begin{lemma}\label{202308091510}
Let $m\geq3$. Then 
\beas
E\big[\Gamma^{(m)}(\sfF_1,...,\sfF_m)\big]
&=&
O(T^{-\frac{m}{2}})1_{\{H\in(\half,\frac{m+1}{2m})\}}
+
O(T^{-\frac{m}{2}+})1_{\{H=\frac{m+1}{2m}\}}
\nn\\&&
+
O(T^{-\frac{m}{2}(3-4H)})1_{\{H\in(\frac{m+1}{2m},1)\}}
\eeas 
as $T\to\infty$, for any $(\sfF_1,...,\sfF_m)\in\{U_T,V_T\}^m$, if 
$\#\{i\in\{1,...,m\};\>\sfF_i=V_T\}=1$. 
\end{lemma}
\proof 
Suppose that $m\geq3$ and $\#\{i\in\{1,...,m\};\>\sfF_i=V_T\}=1$. 
Then we have 
\bea\label{202308091911}
E\big[\Gamma^{(m)}(\sfF_1,...,\sfF_m)\big]
&=& 
2^{m-1}\langle u_T\otimes_1\cdots\otimes_1u_T,v_T\rangle_{\mfh^{\otimes2}}
\yeq
C(m,\theta,H)T^{-m/2}J_T^*,
\eea
where $C(m,\theta,H)$ is a constant and 
\beas 
J_T ^*
&=& 
\int_{[0,T]^{2m}}e^{-\theta|t_1-s_1|}|s_1-t_2|^{2H-2}
e^{-\theta|t_2-s_2|}|s_2-t_3|^{2H-2}
\nn\\&&
\cdots e^{-\theta|t_{m-1}-s_{m-1}|}|s_{m-1}-t_m|^{2H-2}
e^{-\theta|T-t_m|-\theta|T-s_m|}|s_m-t_1|^{2H-2}ds_1dt_2\cdots ds_mdt_1. 
\eeas
\halflineskip
\noindent
1) Case $H\in(\half,\frac{m+1}{2m})$. 
By using (\ref{202307060939}), and (\ref{202307050606}) of Lemma \ref{202307060258}, 
we obtain 
\bea\label{202308091836}
J_T^* 
&\simleq&
\int_{[0,T]^{m+1}}\big(1\wedge|t_1-t_2|^{2H-2}\big)\big(1\wedge|t_2-t_3|^{2H-2}\big)
\cdots\big(1\wedge|t_{m-2}-t_{m-1}|^{2H-2}\big)
\nn\\&&\hspace{30pt}\times
\big(1\wedge|t_{m-1}-t_m|^{2H-2}\big)
e^{-\theta|T-t_m|-\theta|T-s_m|}|s_m-t_1|^{2H-2}dt_1\cdots dt_mds_m
\nn\\&\simleq&
\int_{[0,T]^{m-1}}\big(1\wedge|t_1-t_2|^{2H-2}\big)\big(1\wedge|t_2-t_3|^{2H-2}\big)
\cdots\big(1\wedge|t_{m-2}-t_{m-1}|^{2H-2}\big)
\nn\\&&\hspace{30pt}\times
\big(1\wedge|t_{m-1}-T|^{2H-2}\big)
\big(1\wedge|T-t_1|^{2H-2}\big)dt_1\cdots dt_{m-1}
\nn\\&=&
\int_{[0,T]^{m-1}}\big(1\wedge|t_1-t_2|^{2H-2}\big)\big(1\wedge|t_2-t_3|^{2H-2}\big)
\cdots\big(1\wedge|t_{m-2}-t_{m-1}|^{2H-2}\big)
\nn\\&&\hspace{30pt}\times
\big(1\wedge t_{m-1}^{2H-2}\big)
\big(1\wedge t_1^{2H-2}\big)dt_1\cdots dt_{m-1}. 
\eea
We will estimate the right-hand side of (\ref{202308091836}).
By the same reasoning as the proof of $I_\infty'<\infty$ around (\ref{202308091838}) 
by Young's inequality and H\"older's inequality.   
we see $J_T^*=O(1)$. 
Hence $E\big[\Gamma^{(m)}(\sfF_1,...,\sfF_m)\big]=O(T^{-m/2})$. \y
\noindent
2) Case $H=\frac{m+1}{2m}$. 
For an estimation of the right-hand side of (\ref{202308091836}), 
we can follow the proof of $\wt{I}'_\infty<\infty$ around (\ref{202308091846}), with a discounted function $\wt{a}$. 
Therefore we obtain $E\big[\Gamma^{(m)}(\sfF_1,...,\sfF_m)\big]=O(T^{-\frac{m}{2}+})$. \y
\noindent
3) Case $H\in(\frac{m+1}{2m},1)$. 
Since $|T-t_m|+|T-s_m|\geq|t_m-s_m|$, we have 
\beas 
J_T ^*
&\leq& 
\int_{[0,T]^{2m}}e^{-\theta|t_1-s_1|}|s_1-t_2|^{2H-2}
e^{-\theta|t_2-s_2|}|s_2-t_3|^{2H-2}
\nn\\&&
\cdots e^{-\theta|t_{m-1}-s_{m-1}|}|s_{m-1}-t_m|^{2H-2}
e^{-\theta|t_m-s_m|}|s_m-t_1|^{2H-2}ds_1dt_2\cdots ds_mdt_1
\nn\\&=&
\int_{[0,T]^m}a_T(t_1,t_2)a_T(t_2,t_3)\cdots a_T(t_m,t_1)dt_1\cdots dt_m, 
\eeas
where the function $a_T$ is defined in (\ref{202308091907}). 
Now Lemma \ref{202308120342} gives the estimate 
$J_T ^*=O(T^{m(2H-1)})$, and hence 
$E\big[\Gamma^{(m)}(\sfF_1,...,\sfF_m)\big]
=O(T^{m(2H-3/2)})$ from (\ref{202308091911}). 
This completes the proof of Lemma \ref{202308091510}
\qed\halflineskip

\begin{lemma}\label{202308110241}
Let $H\in(1/2,3/4)$. 
Suppose that $m\geq2$ and $1\leq k\leq m$. Then 
\beas
E\big[\Gamma^{(m)}(\sfF_1,...,\sfF_m)\big]
&=&
O(T^{-\frac{k}{2}})
\eeas 
as $T\to\infty$, for any $(\sfF_1,...,\sfF_m)\in\{U_T,V_T\}^m$, if 
$\#\{i\in\{1,...,m\};\>\sfF_i=V_T\}=k$. 
\end{lemma}
\proof We obtain these estimates from Lemmas \ref{20230813021} and \ref{202308110225}, 
if hypercontractivity and Lemma 3.1 of Tudor and Yoshida \cite{tudor2023high}. 
\qed\halflineskip

\begin{lemma}\label{202308110336}
\bd\im[(a)] $\|W_T\|_{s,p}=O(T^{-1/2})$ as $T\to\infty$ for $s\in\bbR$ and $p>1$. 
\im[(b)] Let $m\geq2$. Then 
$
E\big[\Gamma^{(k)}(\sfF_1,...,\sfF_m)\big]
\yeq
0$
for any $(\sfF_1,...,\sfF_m)\in\{U_T,V_T,W_T\}^m$, if 
$\#\{i\in\{1,...,m\};\>\sfF_i=W_T\}=1$. 
\im[(c)] 
Let $m\geq2$ and $k\leq m$. Then 
$
E\big[\Gamma^{(m)}(\sfF_1,...,\sfF_m)\big]
\yeq
O(T^{-\frac{k}{2}})
$
as $T\to\infty$, for any $(\sfF_1,...,\sfF_m)\in\{U_T,V_T,W_T\}^m$, if 
$\#\{i\in\{1,...,m\};\>\sfF_i=W_T\}=k$. 
\ed
\end{lemma}
\proof (a) is nothing but Lemma \ref{20230813034}. 
(b) follows from the fact that 
$E\big[\Gamma^{(k)}(\sfF_1,...,\sfF_m)\big]$ is the expectation of an element of the first chaos. 
(a) implies (c). 
\qed

\section{Gamma factors and asymptotic expansion of 
{\xred the sum of the basic variables} 
}\label{202402261422}
Define $\bbS_T$ by 
\bea\label{202308090234}
\bbS_T &=& U_T+V_T+W_T, 
\eea
and $c_0$ and $c_2$ by 
\bea\label{202308141446}
c_0\yeq C_U(2,H,\theta)
\quad\text{and}\quad
c_2\yeq {\vred C_U''(2,H,\theta)},
\eea
respectively. 
See (\ref{202308150202}) and (\ref{202308150222}) for these constants. 

\begin{lemma}\label{202308110401} 
Let $H\in(\half,\frac{3}{4})$. Then 
\beas 
E\big[\Gamma^{(2)}(\bbS_T,\bbS_T)\big]
&=& 
c_0+c_2T^{4H-3}+o(T^{4H-3})
\eeas
as $T\to\infty$. 
\end{lemma}
\proof 
From (\ref{202308090234}) and 
{\xred
Lemmas \ref{202308091510a}, 
\ref{202308110241} 
and \ref{202308110336}, 
}
we see
\beas 
E\big[\Gamma^{(2)}(\bbS_T,\bbS_T)\big]
&=& 
E\big[\Gamma^{(2)}(U_T,U_T)\big]+O(T^{-1})
\eeas
as $T\to\infty$. 
We obtain the result from Lemma \ref{202307081346}. 
\qed\halflineskip

Let 
\bea\label{202308141444}
c_3' &=& {\vred C_U(3,H,\theta)}.
\eea
\begin{lemma}\label{202308110413} 
\bd
\im[(a)] For $H\in(\half,\frac{2}{3})$, 
$
E\big[\Gamma^{(3)}(\bbS_T,\bbS_T,\bbS_T)\big]
\yeq
c_3'T^{-\half}
+o\big(T^{-\half}\big). 
$
\im[(b)] 
For $H=\frac{2}{3}$, 
$
E\big[\Gamma^{(3)}(\bbS_T,\bbS_T,\bbS_T)\big]
\yeq 
O(T^{-\half+}). 
$
\im[(c)] For $H\in(\frac{2}{3},1)$, 
$
E\big[\Gamma^{(3)}(\bbS_T,\bbS_T,\bbS_T)\big]
\yeq
O(T^{-\frac{3}{2}(3-4H)}). 
$
\ed
\end{lemma}
\proof 
By using Lemmas \ref{202308091510}, \ref{202308110241} and \ref{202308110336}, 
we obtain 
\beas 
E\big[\Gamma^{(3)}(\bbS_T,\bbS_T,\bbS_T)\big]
&=& 
E\big[\Gamma^{(3)}(U_T,U_T,U_T)\big]
\nn\\&&
+O(T^{-\frac{3}{2}})1_{\{H\in(\half,\frac{2}{3}\}}
+
O(T^{-\frac{3}{2}+})1_{\{H=\frac{2}{3}\}}
\nn\\&&
+
O(T^{-\frac{3}{2}(3-4H)})1_{\{H\in(\frac{2}{3},1)\}}
+
O(T^{-1})
\eeas
as $T\to\infty$. 
Then the desired estimates follow from Lemmas \ref{202307061023}, \ref{202307080000} and \ref{202307061046}.
\qed\halflineskip

The centered $\Gamma^{(p)}$ is denoted by $\wt{\Gamma^{(p)}}$. 
Let 
\beas 
\bbI_T
\yeq
\wt{\Gamma^{(3)}}(\bbS_T,\bbS_T,\bbS_T)
\yeq
\Gamma^{(3)}(\bbS_T,\bbS_T,\bbS_T)-E\big[\Gamma^{(3)}(\bbS_T,\bbS_T,\bbS_T)\big]. 
\eeas

\begin{lemma}\label{202308110455} As $T\to\infty$, 
\beas 
\bbI_T
&=&
1_{\{H\in(\half,\frac{7}{12})\}}O_{D_\infty}(T^{-1})
+1_{\{H=\frac{7}{12}\}}O_{D_\infty}(T^{-1+})
\nn\\&&
+1_{\{H\in(\frac{7}{12},1)\}}O_{D_\infty}(T^{\frac{3}{2}(4H-3)}). 
\eeas
\end{lemma}
\proof 
(I) Estimation of the centered third-order gamma factors 
involving $U_T$ and $V_T$. 
It holds that 
\bea\label{202308130215}
E\big[\big(\wt{\Gamma^{(3)}}(U_T,U_T,U_T)\big)^2\big] 
&=& 
2^4E\big[I_2(u_T\otimes_1u_T\otimes_1u_T)^2\big]
\yeq
2^5\langle \underbrace{u_T\otimes_1\cdots\otimes_1u_T}_5,u_T\rangle_{\mfh\otimes2}
\nn\\&=&
E\big[\Gamma^{(6)}(U_T,...,U_T)\big]
\nn\\&=&
1_{\{H\in(\half,\frac{7}{12})\}}O(T^{-2})
+1_{\{H=\frac{7}{12}\}}O(T^{-2+})
+1_{\{H\in(\frac{7}{12},1)\}}O(T^{3(4H-3)})
\nn\\&&
\eea
from Lemmas \ref{202307061023}, \ref{202307080000} and \ref{202307061046}. 
These estimates are enhanced to $D_\infty$, that is, 
\bea\label{202308130157}
\wt{\Gamma^{(3)}}(U_T,U_T,U_T)
&=&
1_{\{H\in(\half,\frac{7}{12})\}}O_{D_\infty}(T^{-1})
+1_{\{H=\frac{7}{12}\}}O_{D_\infty}(T^{-1+})
\nn\\&&
+1_{\{H\in(\frac{7}{12},1)\}}O_{D_\infty}(T^{\frac{3}{2}(4H-3)}). 
\eea
%

For a mixed centered third-order Gamma factor of $U_T$ and $V_T$, we have 
\beas&&
E\big[\big(\wt{\Gamma^{(3)}}(U_T,U_T,V_T)\big)^2\big] 
\nn\\&=& 
2^4E\big[I_2(u_T\otimes_1u_T\otimes_1v_T)^2\big]\quad(\text{tensor symmtrized})
\nn\\&\sim&
\langle \underbrace{v_T\otimes_1u_T\otimes_1\cdots\otimes_1u_T}_5,v_T\rangle_{\mfh^{\otimes2}}
+\cdots+\langle \underbrace{u_T\otimes_1u_T\otimes_1\cdots\otimes_1v_T}_5,v_T\rangle_{\mfh^{\otimes2}}
\nn\\&\simleq&
T^{-3}\int_{[0,T]^{12}}a(t_1,s_1,t_2)a(t_2,s_2,t_3)\cdots a(t_5,s_5,t_6)a(t_6,s_6,t_1)
dt_1\cdots dt_6ds_1\cdots ds_6. 
\eeas
Here we used $|T-x|+|T-y|\geq|x-y|$ for one $v_T$ to alter it into the function $a$. 
Since 
\beas 
E\big[\big(\wt{\Gamma^{(3)}}(U_T,U_T,V_T)\big)^2\big] 
&\simleq&
E\big[\wt{\Gamma^{(6)}}(U_T,...,U_T)\big] 
\eeas
by (\ref{202307080224}) and (\ref{202307061020}), $\wt{\Gamma^{(3)}}(U_T,U_T,V_T)$ 
admits the same estimate as (\ref{202308130215}), and hence the estimate (\ref{202308130157}). 
\begin{en-text}
Consequently, 
\bea\label{202308130218}
\wt{\Gamma^{(3)}}(V_T,U_T,U_T)
&=&
\wt{\Gamma^{(3)}}(U_T,V_T,U_T)
\yeq
\wt{\Gamma^{(3)}}(U_T,U_T,V_T)
\nn\\&=&
1_{\{H\in(\half,\frac{7}{12})\}}O_{D_\infty}(T^{-1})
+1_{\{H=\frac{7}{12}\}}O_{D_\infty}(T^{-1+})
\nn\\&&
+1_{\{H\in(\frac{7}{12},1)\}}O_{D_\infty}(T^{\frac{3}{2}(4H-3)}). 
\eea
\end{en-text}
On the other hand, Lemmas \ref{20230813021} and \ref{202308110225} give 
$\wt{\Gamma^{(3)}}(V_T,V_T,V_T)=O_{D_\infty}(T^{-3/2})$ and 
$\wt{\Gamma^{(3)}}(V_T,V_T,U_T)=
\wt{\Gamma^{(3)}}(V_T,U_T,V_T)=\wt{\Gamma^{(3)}}(U_T,V_T,V_T)=O_{D_\infty}(T^{-1})$. 
In conclusion, 
\bea\label{202308130227}
\wt{\Gamma^{(3)}}(U'_T,U''_T,U'''_T)
&=&
1_{\{H\in(\half,\frac{7}{12})\}}O_{D_\infty}(T^{-1})
+1_{\{H=\frac{7}{12}\}}O_{D_\infty}(T^{-1+})
\nn\\&&
+1_{\{H\in(\frac{7}{12},1)\}}O_{D_\infty}(T^{\frac{3}{2}(4H-3)}). 
\eea
for $U'_T,U''_T,U'''_T\in\{U_T,V_T\}$.

\noindent
(II) Estimation of the centered third-order gamma factors involving at least one $W_T$. 
We consider $\wt{\Gamma^{(3)}}(U'_T,U''_T,W_T)$ for $U'_T,U''_T\in\{U_T,V_T\}$. 
In order to estimate $E\big[\big(\wt{\Gamma^{(3)}}(U'_T,U''_T,W_T)\big)^2\big]$, 
it suffices to show
\bea\label{202308121455}
\langle \underbrace{w_T\otimes_1u_T\otimes_1\cdots\otimes_1u_T}_k,
\underbrace{w_T\otimes_1u_T\otimes_1\cdots\otimes_1u_T}_{6-k}\rangle_{\mfh}
&=& 
O(T^{-3})
\eea
for $k=0,1,...,5$. 
Here we used the domination of the kernel of $v_T$ by that of $u_T$, once again. 
We also notice that $e^{-2\theta T+\theta t}\leq e^{-\theta t}$ for $t\in[0,T]$. 
Therefore, it is sufficient to use the following estimates: 
\bea
J^{**}_T&:=&
T^{-3}\int_{[0,T]^{9}}a(t,s_1,t_1)a(t_1,s_2,t_2)\cdots a(t_{k-1},s_k,t_k)e^{-\theta t_k}
\nn\\&&\hspace{40pt}\times
a(t,s_{k+1},t_{k+1})a(t_{k+1},s_{k+2},t_{k+2})\cdots a(t_3,s_4,t_4)e^{-\theta t_4}
\nn\\&&\hspace{40pt}\times
dt_1\cdots dt_4ds_1\cdots ds_4dt
\nn\\&\simleq&
T^{-3}\int_{[0,T]^3}(1\wedge|r_1|^{2H-2})(1\wedge|r_1-r_2|^{2H-2})(1\wedge|r_2-r_3|^{2H-2})
(1\wedge|r_3|^{2H-2})
\nn\\&&\hspace{40pt}\times
dr_1dr_2dr_3\label{2023081719}
\\&\simleq&
T^{-(2-\ep)}\int_{[0,T]^3}(1\wedge|r_1|^{2H_1-2})(1\wedge|r_1-r_2|^{2H_1-2})(1\wedge|r_2-r_3|^{2H_1-2})
(1\wedge|r_3|^{2H_1-2})
\nn\\&&\hspace{40pt}\times
dr_1dr_2dr_3,\label{2023081656}
\eea
where $H_1=H-(1+\ep)/8$, $\ep\geq-1$, and $T\geq1$. The last inequality of (\ref{2023081656}) is verified by the estimate
\beas 
T^{-\frac{1+\ep}{4}}(1\wedge|r|^{2H-2})
&\leq&
1\wedge\big(|r|^{-\frac{1+\ep}{4}}|r|^{2H-2}\big)
\yeq
1\wedge|r|^{2H_1-2}
\eeas
for $r\in[-T,T]\setminus\{0\}$ and $T\geq1$. 

When $H\in(\frac{5}{8},\frac{3}{4})$, take $\ep=$ to have $H_1\in(\half,\frac{5}{8})$. 
We apply Lemma \ref{202308120413} to $\alpha_i(x)=1\wedge|x|^{2H_1-2}$ in the case 
$m=4$ and $H_1$ for $H$ under $\ep=0$, to verify the integral on the right-hand side of (\ref{2023081656}) is finite. 
Hence $J^{**}_T=O(T^{-2})$. 

When $H=\frac{5}{8}$, it is possible to show that the integral on the right-hand side of (\ref{2023081656}) is finite 
for any $\ep\in(-1,\infty)$. Therefore, $J^{**}_T=O(T^{-3+})$. 

When $H\in(\half,\frac{5}{8})$, we directly apply Lemma \ref{202308120413} to $\alpha_i(x)=1\wedge|x|^{2H-2}$ in the case 
$m=4$ and $H$, and see integral on the right-hand side of (\ref{2023081719}) is finite, therefore, 
$J^{**}_T=O(T^{-3})$. 

Consequently, for any $H\in(\half,\frac{3}{4})$, $J^{**}_T=O(T^{-2})$, which implies 
$\wt{\Gamma^{(3)}}(U'_T,U''_T,W_T)=O_{D_\infty}(T^{-1})$ as $T\to\infty$, for $U'_T,U''_T\in\{U_T,V_T\}$. 
In the same fashion, it is possible to show $\wt{\Gamma^{(3)}}(U'_T,W_T,U''_T)=O_{D_\infty}(T^{-1})$ and 
$\wt{\Gamma^{(3)}}(W_T,U'_T,U''_T)=O_{D_\infty}(T^{-1})$ for $U'_T,U''_T\in\{U_T,V_T\}$. 

Moreover, Lemmas \ref{20230813021}-\ref{20230813034} show 
$\wt{\Gamma^{(3)}}(W_T,W_T,U'_T)$, $\wt{\Gamma^{(3)}}(W_T,U'_T,W_T)$ and $\wt{\Gamma^{(3)}}(U'_T,W_T,W_T)$ 
are of order $O_{D_\infty}(T^{-1})$ for $U'_T\in\{U_T,V_T\}$. 
Similarly, $\wt{\Gamma^{(3)}}(W_T,W_T,W_T)=O_{D_\infty}(T^{-3/2})$. 

After all that, 
\bea\label{202308130245}
\wt{\Gamma^{(3)}}(U'_T,U''_T,U'''_T)
&=& 
O_{D_\infty}(T^{-1})
\eea
for $U'_T,U''_T\in\{U_T,V_T,W_T\}$ if $1_{\{U'_T=W_T\}}+1_{\{U''_T=W_T\}}+1_{\{U'''_T=W_T\}}\geq1$. 

\noindent
(III) The proof of Lemma \ref{202308110455} is completed by merging 
(\ref{202308130227}) and (\ref{202308130245}). 
\qed\halflineskip

The estimated exponents of $T$ and the ranks of the terms appearing in the asymptotic expansion are summarized 
in Table \ref{202308141322}, together with the estimates for the centered third-order gamma factors. 
It should be remarked that the change of the second dominant terms is seamless at $H=5/8$. 
In the asymptotic expansion, the classical order $-1/2$ becomes the exponent of the first-order correction term for $H\in(1/2,5/8)$, 
while $4H-3$ does for $H\in(5/8,3/4)$, and both do at $H=5/8$. 
\begin{table}[H]
  \caption{Estimated exponents of $T$ and [Rank]s}
  \label{202308141322}
  \centering
  \begin{tabular}{lcccc}
    \hline
    sequence\textbackslash interval  & $(\frac{1}{2},\frac{5}{8})$ &$(\frac{5}{8},\frac{7}{12})$ &$(\frac{7}{12},\frac{2}{3})$ &   $(\frac{2}{3},\frac{3}{4})$  \\
    \hline \hline
    0th-order term of $E[\Gamma^{(2)}(U_T,U_T)]$  & 0\>[1] &0\>[1]&0\>[1]& 0\>[1] \\
    1st-order term of $E[\Gamma^{(2)}(U_T,U_T)]$   & $4H-3\>[3]$  &$4H-3\>[2]$&$4H-3\>[2]$& $4H-3\>[2]$ \\
     $E[\Gamma^{(3)}(\bbS_T,\bbS_T,\bbS_T)]$  &$-\half\>[2]$& $-\half\>[3]$  &$-\half\>[3]$& $\frac{3}{2}(4H-3)\>[3]$ \vspace{3pt}\\
     \hline
    $E[\wt{\Gamma}^{(3)}(U_T,U_T,U_T)]$  &-1& $-1$  &$\frac{3}{2}(4H-3)$& $\frac{3}{2}(4H-3)$ \\
    $E[\wt{\Gamma}^{(3)}(U_T,U_T,V_T)]$  &-1& $-1$  &$\frac{3}{2}(4H-3)$& $\frac{3}{2}(4H-3)$ \\
     $E[\wt{\Gamma}^{(3)}(U_T',U_T'',W_T)]$  &-1& $-1$  &$-1$& $-1$ \\    \hline
  \end{tabular}
\end{table}

\begin{comment}
Comment *********************\\
The order  $\sfq$ of the expansion is defined as 
\beas
\sfq\yeq\sfq(H)
&=&
\l\{\begin{array}{ll}
\half&\big(H\in\big(\frac{1}{2},\frac{5}{8}\big]\big)\y
-4H+3&\big(H\in\big(\frac{5}{8},\frac{3}{4}\big)\big)
\end{array}\r.
\eeas
**********************************
\end{comment}

We shall derive an asymptotic expansion of $\bbS_T$. 
Define the density function $p^*_{H,T}(x)$ as 
\bea\label{202308141711}
p^*_{H,T}(x)
&=& 
\phi(x;0,c_0)\bigg(1
+1_{\{H\in[\frac{5}{8},\frac{3}{4})\}}2^{-1}c_2H_2(x;0,c_0)T^{4H-3}
\nn\\&&\hspace{63pt}
+1_{\{H\in(\half,\frac{5}{8}]\}}3^{-1}c_3'H_3(x;0,c_0)T^{-\half}\bigg).
\eea
The exponent $\sfq=\sfq(H)$ is given in (\ref{202307110532}). 
\begin{proposition}\label{202308141705}
Suppose that $H\in(1/2,3/4)$. Then 
\bea\label{202308141640}
\sup_{g\in\cale(a,b)}\bigg|E[g(\bbS_T)-\int_\bbR g(x)p^*_{H,T}(x)dx\bigg|
&=& 
o(T^{-\sfq(H)})
\eea
as $T\to\infty$. 
\end{proposition}
\proof 
Prepare the following parameters: 
\beas 
d\yeq1,\quad
p\yeq2,\quad
\sfk \yeq 1, \quad
\sfq_0(H) \yeq \frac{2}{3}\sfq(H), \quad \xi(H)\yeq \frac{1}{9}\sfq(H), \quad 
\ell\yeq11,\quad\ell_1\yeq5. 
\eeas
Then 
\beas &&
\sfq_0(H)(\sfk+1)>\sfq(H),\quad
\xi(H)(\ell-d)\>>\>\sfq(H),\quad
\nn\\&&
\ell\geq\ell_1>p+1+d,\quad
\sfq_0(H) \yleq \sfq(H)-3\xi(H). 
\eeas
Therefore, Condition $[B]$ of Tudor and Yoshida \cite{tudor2023high} is satisfied for each $H\in(1/2,3/4)$, 
thanks to Lemmas \ref{202308110401} and \ref{202308110413}. 

\begin{en-text}
\beas 
\bbS_T
&=& 
{\bf G}(\theta)^{-1}S_T
-T^{-1/2}{\bf G}(\theta)^{-3}{\bf C}(\theta)S_T^2. 
\eeas
\end{en-text}

From (\ref{202307030144}), 
the formula (\ref{202307030127}) gives 
\bea\label{202308141618}
\Gamma^{(2)}(U_T,U_T) 
&=& 
2I_2\big(u_T\otimes_1u_T\big)
+2\big\langle u_T,u_T\big\rangle_{\calh^{\otimes2}}. 
\eea
Lemma \ref{202307081346} shows 
\bea\label{202308141619}
2\big\langle u_T,u_T\big\rangle_{\calh^{\otimes2}}
&=&
C_U(2,H,\theta)+O(T^{4H-3}). 
\eea
From (\ref{202308141618}) and (\ref{202308141619}), 
\beas
\Gamma^{(2)}(U_T,U_T) - c_0 
&=& 
2I_2\big(u_T\otimes_1u_T\big)+O(T^{4H-3})
\eeas
Furthermore, 
\beas 
E\big[I_2\big(u_T\otimes_1u_T\big)^2\big]
&=&
2\langle u_T\otimes_1u_T,u_T\otimes_1u_T\rangle_{\mfh^{\otimes2}}
\nn\\&=&
2\langle u_T\otimes_1u_T\otimes_1u_T,u_T\rangle_{\mfh^{\otimes2}}
\nn\\&=&
1_{\{H\in(\half,\frac{5}{8}\}}O(T^{-1})+1_{\{H=\frac{5}{8}\}}O(T^{-1+})+1_{\{H\in(\frac{5}{8},\frac{3}{4})\}}O(T^{2(4H-3)}). 
\eeas
by Lemmas \ref{202307061023}, \ref{202307080000} and \ref{202307061046}. 
Therefore, in any case of $H\in(1/2,3/4)$, we can find a positive constant $\sfa(H)$ such that 
\beas
\Gamma^{(2)}(U_T,U_T) - c_0 
&=& 
O_{D_\infty}(T^{-\sfa(H)})
\eeas
as $T\to\infty$. 
With the help of Lemmas \ref{202308110225} and \ref{20230813034}, this verifies $[A1]$ (ii) of Tudor and Yoshida \cite{tudor2023high} for $\Gamma^{(2)}(\bbS_T,\bbS_T)$. 
Lemmas \ref{20230813021}-\ref{20230813034} imply $\bbS_T=O_{D_\infty}(1)$, and $[A1]$ (i) is checked. 
Thus, $[A1]$ of Tudor and Yoshida \cite{tudor2023high} holds. 
Besides, Condition $[A2^\sharp]$ of Tudor and Yoshida \cite{tudor2023high} has been ensured by Lemma \ref{202308110455}. 
We apply Theorem 5.2 of Tudor and Yoshida \cite{tudor2023high} to conclude (\ref{202308141640}). 
\qed

\section{Smooth stochastic expansion of the estimator}\label{202402261439}
Let $Q_T=\int_0^TX_t^2dt$. 
Define ${\bf G}(\vartheta)$ by 
\bea\label{202308140219}
{\bf G}(\vartheta)
&=& 
\int_0^1\partial_\theta\mu\big(\theta+u(\vartheta-\theta)\big)du\quad(\vartheta\in(0,\infty)). 
\eea
In particular, 
\bea\label{202308140222}
{\bf G}(\theta)
\yeq
\partial_\theta\mu(\theta)
\yeq
-2\sigma^2H^2\Gamma(2H)\theta^{-2H-1}.
\eea
\begin{lemma}\label{202308131432}
\bea\label{202308131614}
Q_T
&=& 
T^{1/2}{\bf G}(\theta)\big(U_T+V_T+W_T)+\ol{\nu}_T(\theta). 
\eea
\end{lemma}
\proof
By the representation 
\beas 
X_t &=& e^{-\theta t}x_0+I_1\big(\sigma e^{-\theta(t-\cdot)}1_{[0,t]}(\cdot)\big), 
\eeas
we have
\bea\label{202308131449} 
X_t^2 
&=& 
e^{-2\theta t}x_0^2+2e^{-\theta t}x_0I_1\big(\sigma e^{-\theta(t-\cdot)}1_{[0,t]}(\cdot)\big)
\nn\\&&
+I_2\bigg(\sigma^2e^{-\theta(t-\cdot)}1_{[0,t]}(\cdot)\otimes e^{-\theta(t-\cdot)}1_{[0,t]}(\cdot)\bigg)
+\sigma^2\big\langle e^{-\theta(t-\cdot)}1_{[0,t]}(\cdot), e^{-\theta(t-\cdot)}1_{[0,t]}(\cdot)\big\rangle_\mfh
\nn\\&=&
e^{-2\theta t}x_0^2+2e^{-\theta t}x_0I_1\big(\sigma e^{-\theta(t-\cdot)}1_{[0,t]}(\cdot)\big)
+I_2\bigg(\sigma^2e^{-\theta(t-\cdot)}1_{[0,t]}(\cdot)\otimes e^{-\theta(t-\cdot)}1_{[0,t]}(\cdot)\bigg)
\nn\\&&
+\sigma^2\alpha_H
\int_{[0,t]^2}
e^{-\theta(t-s_1)}e^{-\theta(t-s_2)}|s_1-s_2|^{2H-2}ds_1ds_2. 
\eea
Moreover, 
\bea\label{202308131518}&&
\int_0^T \sigma^2e^{-\theta(t-s_1)}1_{[0,t]}(s_1)e^{-\theta(t-s_2)}1_{[0,t]}(s_2)dt
\nn\\&=&
\int_{s_1\vee s_2}^T\sigma^2 e^{-2\theta t+\theta(s_1+s_2)}dt1_{\{s_1, s_2\in[0,T]\}}
\nn\\&=&
\sigma^2(2\theta)^{-1} \big(
e^{-\theta |s_1-s_2|}-e^{\theta(-2 T+s_1+s_2)}\big)1_{\{s_1, s_2\in[0,T]\}}
\nn\\&=&
T^{1/2}\sigma^2(2\theta K_U)^{-1}u_T(s_1,s_2)-T^{1/2}\sigma^2(2\theta K_V)^{-1}v_T(s_1,s_2),
\eea
and
\bea\label{202308131519}
\int_0^T 2x_0\sigma e^{-2\theta t+\theta s}1_{\{s<t\leq T\}}dt
&=&
x_0\sigma\theta^{-1}\big(e^{-\theta s}-e^{-2\theta T+\theta s}\big)1_{\{s\in[0,T]\}}
\nn\\&=&
T^{1/2}x_0\sigma\theta^{-1}K_W^{-1}w_T(s). 
\eea
Therefore, 
(\ref{202308131449}), (\ref{202308131518}) and (\ref{202308131519}) gives (\ref{202308131614}). 
\qed\halflineskip

\begin{comment}
Comment\hrulefill
\beas&&
\int_0^T \sigma^2e^{-\theta(t-s_1)}1_{[0,t]}(s_1)e^{-\theta(t-s_2)}1_{[0,t]}(s_2)dt
\nn\\&=&
\int_{s_1\vee s_2}^T\sigma^2 e^{-2\theta t+\theta(s_1+s_2)}dt1_{\{s_1, s_2\in[0,T]\}}
\nn\\&=&
\sigma^2(2\theta)^{-1} \big(
e^{-2\theta |s_1-s_2|}-e^{\theta(-2 T+s_1+s_2)}\big)1_{\{s_1, s_2\in[0,T]\}}
\nn\\&=&
T^{1/2}\theta^{-1}f_T(s_1,s_2)-h_T(s_1,s_2)
\eeas
\hrulefill
\end{comment}

%
\begin{lemma}\label{202308131351}
For every $\ep>0$ and $L>0$, 
$P\big[|\wh{\theta}_T-\theta|>\ep\big]=O(T^{-L})$ as $T\to\infty$. 
\end{lemma}
\proof 
Take a sufficiently small positive number $\sfr$ such that $U(\theta,\sfr)\equiv\{\theta'\in\bbR;|\theta'-\theta|<\sfr\}\subset\Theta$. Suppose that $0<2\ep<\sfr$. 
By definition of $\wh{\theta}_T$, we have 
\bea\label{202308131723}
\big\{|\wh{\theta}_T-\theta|>2\ep\big\}
&\subset&
\big\{\big|\wt{\theta}_T-\theta\big|>\ep\big\}\cup\big\{T^{-1}\|\beta\|_\infty>\ep\big\}%
\nn\\&\subset&
\bigg\{|T^{-1}Q_T-\mu(\theta)|\geq\inf_{\theta':|\theta'-\theta|>\ep}|\mu(\theta')-\mu(\theta)|\bigg\}
\cup\big\{T^{-1}\|\beta\|_\infty>\ep\big\}%
\eea
since $T^{-1}Q_T=\mu(\wt{\theta}_T)$. 
Then 
\beas
P\big[|\wh{\theta}_T-\theta|>\ep\big]
&\simleq&
E\big[\big|T^{-1}Q_T-T^{-1}\ol{\nu}_T(\theta)\big|^{2L}\big]
\yeq
O(T^{-L})
\eeas
as $T\to\infty$ {\ared(recall $\ol{\nu}_T(\theta)=E\big[\int_0^TX_t^2dt\big]$)} 
since $T^{-1/2}(Q_T-\ol{\nu}_T(\theta))=O_{L^\inftym}(1)$ as $T\to\infty$, i.e., 
all $L^p$-norms are bounded, from 
the representation (\ref{202308131614}) of $Q_T$ and Lemmas \ref{20230813021}-\ref{20230813034}. 
\qed\halflineskip

\begin{comment}
\hrulefill\\
Comment 
By definition of $\wh{\theta}_T$, we have 
\beas
\big\{[|\wh{\theta}_T-\theta|>\ep\big\}
&\subset&
\bigg\{\inf_{\theta':|\theta'-\theta|>\ep}\big|\psi_T(\theta')\big|\leq\big|\psi_T(\theta)\big|\bigg\}
\nn\\&\leq&
\bigg[\{\inf_{\theta':|\theta'-\theta|>\ep}\big|T^{-1}Q_T-T^{-1}\nu_T(\theta')\big|\leq\big|T^{-1}Q_T-T^{-1}\nu_T(\theta)\big|\bigg\}
\nn\\&\subset&
\bigg\{[\inf_{\theta':|\theta'-\theta|>\ep}\big|T^{-1}\nu_T(\theta')-T^{-1}\ol{\nu}_T(\theta)\big|
\nn\\&&\hspace{40pt}\leq
2\big|T^{-1}Q_T-T^{-1}\ol{\nu}_T(\theta)\big|+\big|T^{-1}\ol{\nu}_T(\theta)-T^{-1}\nu_T(\theta)\big|
\bigg\}
\nn\\&\subset&
\bigg\{\inf_{\theta':|\theta'-\theta|>\ep}\big|T^{-1}\wt{\nu}_T(\theta')-T^{-1}\wt{\nu}_T(\theta)\big|
-2T^{-1}|\ol{b}_T(\theta)|-2T^{-1}\sup_{\theta'\in\Theta}|b_T(\theta')|
\nn\\&&\hspace{40pt}\leq
2\big|T^{-1}Q_T-T^{-1}\ol{\nu}_T(\theta)\big|
\bigg\}.
\eeas
Then 
\beas
P\big[|\wh{\theta}_T-\theta|>\ep\big]
&\simleq&
E\big[\big|T^{-1}Q_T-T^{-1}\ol{\nu}_T(\theta)\big|^{2L}\big]
\eeas
since $|\ol{b}_T(\theta)|\koko+\sup_{\theta'\in\Theta}|b_T(\theta')|<\infty$, where 
$\ol{b}_T(\theta)=\ol{\nu}_T(\theta)-\wt{\nu}_T(\theta)$. 
%
Furthermore, we obtain $T^{-1/2}(Q_T-\ol{\nu}_T(\theta))=O_{L^\inftym}(1)$ as $T\to\infty$, i.e., 
all $L^p$-norms are bounded, from 
the representation (\ref{202308131614}) of $Q_T$ and Lemmas \ref{20230813021}-\ref{20230813034}. 
This proves the result.
\\ \hrulefill\\
\end{comment}

Let 
\bea\label{202308141909} 
b_\infty(\theta)
&=& 
-\sigma^2\alpha_H\Gamma(2H)\theta^{-2H-1}
{\bred -\half\sigma^2\alpha_H\Gamma(2H-1)\theta^{-2H-1}}
+\frac{1}{2\theta}x_0^2
\nn\\&=&
{\bred  -\half\sigma^2\alpha_H(4H-1)\Gamma(2H-1)\theta^{-2H-1}+\frac{1}{2\theta}x_0^2}
\eea

\begin{lemma}\label{202308132345}
$\ol{\nu}_T(\theta)=\wt{\nu}_T(\theta)+\ol{b}_T(\theta)$ and $\ol{b}_T(\theta)\to b_\infty(\theta)$ as $T\to\infty$. 
\end{lemma}
\proof We see 
\beas 
\ol{\nu}_T(\theta)
&=&
E\bigg[\int_0^TX_t^2dt\bigg]
\nn\\&=&
2\sigma^2\alpha_H(2\theta)^{-1}\int_0^Te^{-\theta t}t^{2H-2}dt\>T
-2\sigma^2\alpha_H(2\theta)^{-1}\int_0^Tte^{-\theta t}t^{2H-2}dt
\nn\\&&
-\sigma^2\alpha_H\int_{[0,T]^2}(2\theta)^{-1}e^{-\theta (s_1+s_2)}|s_1-s_2|^{2H-2}ds_1ds_2
+\frac{1-e^{-2\theta T}}{2\theta}x_0^2
\nn\\&=&
\wt{\nu}_T(\theta)+\ol{b}_T(\theta).
\eeas
Remark that 
\beas
2\alpha_H(2\theta)^{-1}\int_0^Te^{-\theta t}t^{2H-2}dt
&=&
H(2H-1)\Gamma(2H-1)\theta^{-2H}+O(e^{-\theta T/2})
\nn\\&=&
H\Gamma(2H)\theta^{-2H}+O(e^{-\theta T/2})
\eeas
as $T\to\infty$. 
Therefore, 
\beas 
\lim_{T\to\infty}\ol{b}_T(\theta)
&=& 
-2\sigma^2\alpha_H(2\theta)^{-1}\int_0^\infty te^{-\theta t}t^{2H-2}dt
\nn\\&&
-\sigma^2\alpha_H\int_{[0,\infty)^2}(2\theta)^{-1}e^{-\theta (s_1+s_2)}|s_1-s_2|^{2H-2}ds_1ds_2
+\frac{1}{2\theta}x_0^2
\nn\\&=&
-\sigma^2\alpha_H\Gamma(2H)\theta^{-2H-1}
{\bred -\half\sigma^2\alpha_H\Gamma(2H-1)\theta^{-2H-1}}
+\frac{1}{2\theta}x_0^2.
\eeas
The proof is completed. 
\qed\halflineskip

The effect of the initial value $x_0$ may appear in the asymptotic expansion 
possibly in the leading correction term. 
In this sense, we can say the moment estimator is fairly skewed.

When $\wt{\theta}_T\in U(\theta,\sfr)$ and $\wh{\theta}_T^o\in U(\theta,\sfr)$, 
\bea\label{202308080345}
S_T &:=& 
T^{-1/2}\big(Q_T-\ol{\nu}_T(\theta) \big)
\nn\\&=& 
T^{-1/2}\big(\wt{\nu}_T(\wt{\theta}_T)-\ol{\nu}_T(\theta) \big)
\nn\\&=&
{\bf G}(\wt{\theta}_T)\>T^{1/2}(\wt{\theta}_T-\theta) -T^{-1/2}\ol{b}_T(\theta)
\eea
and 
\bea\label{202308080346}
S_T 
&=& 
{\bf G}(\theta)\>T^{1/2}(\wt{\theta}_T-\theta) 
+T^{-1/2}{\bf C}(\wt{\theta}_T)\>T(\wt{\theta}_T-\theta)^2
-T^{-1/2}\ol{b}_T(\theta),
\eea
where ${\bf G}(\vartheta)$ is defined by (\ref{202308140219}) and 
\beas 
{\bf C}(\vartheta)
&=&
\int_0^1(1-u)\partial_\theta^2\mu\big(\theta+u(\vartheta-\theta)\big)du. 
\eeas
By definition, 
$
{\bf G}(\theta) \yeq 
-2\sigma^2H^2\Gamma(2H)\theta^{-2H-1}
$ (see (\ref{202308140222})) 
and 
\beas 
{\bf C}(\theta) 
\yeq
\sigma^2H^2(2H+1)\Gamma(2H)\theta^{-2H-2}
\yeq
2^{-1}\sigma^2H\Gamma(2H+2)\theta^{-2H-2}. 
\eeas

Since
$
\inf_{\vartheta\in\ol{\Theta}}|{\bf G}(\vartheta)| > 0
$, 
we have 
\bea\label{202308080347}
T^{1/2}(\wt{\theta}_T-\theta) &=& 
{\bf G}(\wt{\theta}_T)^{-1}S_T\> +T^{-1/2}{\bf G}(\wt{\theta}_T)^{-1}\ol{b}_T(\theta)
\eea
from (\ref{202308080345}), besides 
\bea\label{202308080348}
T^{1/2}(\wt{\theta}_T-\theta) &=& 
{\bf G}(\theta)^{-1}S_T
-T^{-1/2}{\bf G}(\theta)^{-1}{\bf C}(\wt{\theta}_T)\>T(\wt{\theta}_T-\theta)^2
\nn\\&&
+T^{-1/2}{\bf G}(\theta)^{-1}\ol{b}_T(\theta)
\eea
from (\ref{202308080346}). 
Substitute the expression of (\ref{202308080347}) for $T(\wt{\theta}_T-\theta)^2$ of (\ref{202308080348}) 
to obtain 
\bea\label{202308081458}
T^{1/2}(\wt{\theta}_T-\theta) 
&=& 
{\bf G}(\theta)^{-1}S_T
-T^{-1/2}{\bf G}(\theta)^{-3}{\bf C}(\theta)S_T^2
\nn\\&&
+T^{-1/2}{\bf G}(\theta)^{-1}\ol{b}_T(\theta)+{\bf R}_T^\dagger,
\eea
where 
\bea\label{202308141802}
{\bf R}_T^\dagger
&=&
-T^{-1/2}{\bf G}(\theta)^{-3}\big({\bf C}(\wt{\theta}_T)-{\bf C}(\theta)\big)S_T^2
\nn\\&&
-T^{-1/2}{\bf G}(\theta)^{-1}{\bf C}(\wt{\theta}_T)\big\{2S_T{\bf R}_T^*+({\bf R}_T^*)^2\big\}
\eea
with ${\bf R}_T^*$ given by 
\bea\label{202308141803}
{\bf R}_T^*(\theta)
&=& 
\big({\bf G}(\wt{\theta}_T)^{-1}-{\bf G}(\theta)^{-1}\big)S_T\> +T^{-1/2}{\bf G}(\wt{\theta}_T)^{-1}\ol{b}_T(\theta).
\eea
Finally, from (\ref{202308081458}), 
\bea\label{202308140412}
{\sred T^{1/2}}(\wh{\theta}_T-\theta)
&=& 
{\bf G}(\theta)^{-1}S_T
-T^{-1/2}{\bf G}(\theta)^{-3}{\bf C}(\theta)S_T^2
\nn\\&&
+T^{-1/2}{\bf G}(\theta)^{-1}\ol{b}_T(\theta)-T^{{\sred-\half-\sfq(H)}}\beta(\theta)+{\bf R}_T^\ddagger, 
\eea
where 
\bea\label{202308141800}
{\bf R}_T^\ddagger
&=& 
{\bf R}_T^\dagger-T^{-1/2}\big(\beta(\wt{\theta}_T)-\beta(\theta)\big). 
\eea

\begin{en-text}

$F_T^{o}$ is explicit of $Q_T$. 

General $F_T$ is a smooth function of $F_T^{o}$, the function is stable in $T$. 

\beas 
\Xi(x,\theta,\beta) 
&=& 
x-\mu(\theta)-\beta^2 b_{1/\beta^2}(\theta)
\eeas
with $\mu(\theta)=\sigma^2H\Gamma(2H)\theta^{-2H}$. 
\end{en-text}
\begin{en-text}
Thus, 
\beas
S_T &=& 
T^{-1/2}\big(Q_T-\ol{\nu}_T(\theta) \big)
\nn\\&=&
x_0\sigma\theta^{-1} I_1\big((e^{-\theta\cdot}-e^{-2\theta T+\theta\cdot})1_{\{\cdot\in[0,T]\}}\big)+I_2\big(\theta^{-1}f_T-T^{-1/2}h_T\big)
\nn\\&=&
x_0\sigma\theta^{-1} I_1\big((e^{-\theta\cdot}-e^{-2\theta T+\theta\cdot})1_{\{\cdot\in[0,T]\}}\big)+\theta^{-1}F_T-T^{-1/2}I_2\big(h_T\big)
\nn\\&=&
x_0\sigma\theta^{-1} I_1\big((e^{-\theta\cdot}-e^{-2\theta T+\theta\cdot})1_{\{\cdot\in[0,T]\}}\big)+T^{1/2}\big(G_T-E[G_T]\big)
\nn\\&=&
x_0\sigma\theta^{-1} I_1\big((e^{-\theta\cdot}-e^{-2\theta T+\theta\cdot})1_{\{\cdot\in[0,T]\}}\big)+T^{1/2}I_2(g_T)
\eeas
\end{en-text}
%

\begin{comment} 
Comment \hrulefill\\
Recall 
\beas 
u_T(s,t) &=&  
{\bf G}(\theta)^{-1}(2\theta)^{-1}\sigma^2T^{-1/2}e^{-\theta|s-t|}1_{[0,T]^2}(s,t) 
\yeq 
K_UT^{-1/2}e^{-\theta|s-t|}1_{[0,T]^2}(s,t)
\nn\\&&
\text{with }K_U \yeq 
-\frac{\theta^{2H}}{4H^2\Gamma(2H)}, 
\nn\\
v_T(s,t) &=& -{\bf G}(\theta)^{-1}(4\theta)^{-1}\sigma^2T^{-1/2}e^{-\theta(T-s)-\theta(T-t)}1_{[0,T]^2}(s,t)
\nn\\&=&
K_VT^{-1/2}e^{-\theta(T-s)-\theta(T-t)}1_{[0,T]^2}(s,t)
\nn\\&&
\text{with }K_V \yeq \frac{\theta^{2H}}{{\vred4}H^2\Gamma(2H)}, 
\nn\\
w_T(t) &=& T^{-1/2}{\bf G}(\theta)^{-1}x_0\sigma\theta^{-1} (e^{-\theta t}-e^{-2\theta T+\theta t})1_{[0,T]}(t)
\yeq
K_WT^{-1/2}(e^{-\theta t}-e^{-2\theta T+\theta t})1_{[0,T]}(t)
\nn\\&&
\text{with }K_W \yeq 
-\frac{x_0\theta^{2H}}{2\sigma H^2\Gamma(2H)}. 
\eeas
%
\beas 
u_T(s,t) &=&  
{\bf G}(\theta)^{-1}(2\theta)^{-1}\sigma^2T^{-1/2}e^{-\theta|s-t|}1_{[0,T]^2}(s,t) (\yeq {\bf G}(\theta)^{-1}\theta^{-1}f_T)
\nn\\
v_T(s,t) &=& -{\bf G}(\theta)^{-1}(4\theta)^{-1}\sigma^2T^{-1/2}e^{-\theta(T-s)-\theta(T-t)}1_{[0,T]^2}(s,t)
\nn\\&&
(\yeq-{\bf G}(\theta)^{-1}2^{-1}T^{-1/2}h_T\text{ or }-{\bf G}(\theta)^{-1}2^{-1}\sigma^2T^{-1/2}h_T)
\nn\\
w_T(t) &=& T^{-1/2}{\bf G}(\theta)^{-1}x_0\sigma\theta^{-1} (e^{-\theta t}-e^{-2\theta T+\theta t})1_{[0,T]}(t). 
\eeas
\hrulefill
\end{comment} 

Recall $\bbS_T={\bf G}(\theta)^{-1}S_T$ has the representation 
\beas
\bbS_T &=& U_T+V_T+W_T. 
\eeas
From (\ref{202308140412}). 
\bea\label{202308090228}
T^{1/2}(\wh{\theta}_T-\theta) 
&=& 
\bbS_T
+T^{-1/2}{\sred\lambda}\bbS_T^2
+T^{-{\sred\sfq(H)}}\bbd_T+{\bf R}_T^\ddagger,
\eea
where 
\beas 
\bbd_T={\sred T^{-\half+\sfq(H)}}{\bf G}(\theta)^{-1}\ol{b}_T(\theta)-\beta(\theta)
\eeas 
and 
\bea\label{202308150236} 
{\sred\lambda}
&=& 
-{\bf G}(\theta)^{-1}{\bf C}(\theta)
\yeq
2^{-1}(2H+1)\theta^{-1}. 
\eea

Take a smooth function $\psi:\bbR\to[0,1]$ such that $\psi(x)=1$ when $|x|<1/2$ and 
$\psi(x)=0$ when $|x|>1$. 
Let 
\bea\label{202308140428}
\psi_T^{C_1}=\psi\big(C_1\big|T^{-1}Q_T-T^{-1}\ol{\nu}_T(\theta)\big|^2\big). 
\eea
In view of (\ref{202308131723}), we can say there exist numbers $T_1$ and $C_1$ such that 
$\wt{\theta}_T\in U(\theta,\sfr)$ and 
$\wh{\theta}_T\in U(\theta,\sfr)$ whenever $\psi_T>0$ and $T>T_1$. 
In what follows, we will only consider $T$ such that $T>T_1$. 
Then the functional $\wt{F}_T^{C_1}:=\psi_T^{C_1}T^{1/2}(\wt{\theta}_T-\theta)$ is well defined on the whole probability space 
and it is possible to show $\wt{F}_T^{C_1}=O_{D_\infty}(1)$. 
In this way, we have reached the stochastic expansion 
\bea\label{202308141731}
F_T^{2C_1}
&:=& 
\psi_T^{2C_1}T^{1/2}(\wh{\theta}_T-\theta) \yeq 
\bbS_T
+T^{-1/2}\kappa\bbS_T^2
+T^{{\sred-\sfq(H)}}\bbd_T+{\bf R}_T,
\eea
where 
\bea\label{202308141806} 
{\bf R}_T
&=& 
\psi_T^{2C_1}{\bf R}_T^\ddagger
-(1-\psi_T^{2C_1})\big(\bbS_T
+T^{-1/2}\kappa\bbS_T^2
+T^{{\sred-\sfq(H)}}\bbd_T\big).
\eea

\begin{lemma}\label{202308141801}
${\bf R}_T\in D_\infty$ and 
${\bf R}_T=O_{D_\infty}(T^{-1})$ as $T\to\infty$. 
\end{lemma}
\proof 
It is easy to show that $\psi_T^{2C_1}\in D_\infty$ and 
$\psi_T^{2C_1}-1=O_{D_\infty}(T^{-L})$ for every $L>0$. 
As for the term $\psi_T^{2C_1}{\bf R}_T^\ddagger$ in (\ref{202308141806}), 
it is observed that, on the event $\{\psi_T^{2C_1}>0\}$, the terms appearing in the representation of ${\bf R}_T^\ddagger$ consist of 
some functionals of the form $f(\wt{\theta}_T)$ for a $f\in C_B^\infty(U(\theta,\sfr))$. 
Since $\psi_T^{2C_1}{\bf R}_T^\ddagger$ has the factor $\psi_T^{2C_1}$, 
we can replace $f(\wt{\theta}_T)$ by $f(\theta+T^{-1/2}\wt{F}^{C_1}_T)$. 
The latter is well defined on the whole probability space and indeed it is in $D_\infty$. 
Along (\ref{202308141800}), (\ref{202308141802}) and (\ref{202308141803}), 
we can verify that ${\bf R}_T\in D_\infty$ and 
${\bf R}_T=O_{D_\infty}(T^{-1})$ as $T\to\infty$. 
\qed\halflineskip

\section{Proof of Theorems \ref{202308141824} and \ref{202403210212}}\label{202402261518}
\subsection{Proof of Theorem \ref{202308141824}}\label{202403210312}
The asymptotic expansion $p^*_{H,T}$ for $\bbS_T$ has already been obtained in Proposition \ref{202308141705}. 
We will deal with 
the last three terms on the right-hand side of (\ref{202308141731}) by the perturbation method 
of Sakamoto and Yoshida \cite{SakamotoYoshida2003}. 
The stochastic expansion (\ref{202308141731}) of $F_T^{2C_1}$ reads 
$F_T^{2C_1}=\bbS_T+{\sred T^{-\sfq(H)}}\bbY_T$ with 
the perturbation term $\bbY_T={\sred T^{\sfq(H)-\half}}\kappa\bbS_T^2+\bbd_T+{\sred T^{\sfq(H)}}{\bf R}_T$. 
{\sblue From Proposition \ref{202308141705}, in particular, }
\beas
(\bbS_T,\bbY_T) &\to^d& \big(\bbS_\infty,
{\sred 1_{\{H\in(\half,\frac{5}{8}]\}}}\kappa\bbS_\infty^2+{\sred 1_{\{H\in(\half,\frac{5}{8}]\}}}{\bf G}(\theta)^{-1}b_\infty(\theta)-\beta(\theta)\big)
\eeas
as $T\to\infty$, where 
{\sred $\bbS_\infty$ is a random variable distributed as $\bbS_\infty\sim N(0,c_0)$ and} 
$b_\infty(\theta)$ is given in (\ref{202308141909}). 
We can apply Theorem 2.1 of Sakamoto and Yoshida \cite{SakamotoYoshida2003} because 
asymptotic non-degeneracy of $\bbS_T$ is obvious. 
The asymptotic expansion for $F_T^{2C_1}$ is now given by the density function 
\bea\label{202308141912}
p_{H,T}(x)
&=& 
p_{H,T}^*(x)+{\sred T^{-\sfq(H)}}g(x),
\eea
where 
\beas 
g(x) 
&=& 
-\partial_x\big\{(\kappa x^2+\tau)\phi(x;0,c_0)\big\}
\eeas
with 
\bea\label{202308150237}
{\sred\kappa\yeq\kappa(H,\theta)\yeq1_{\{H\in(\half,\frac{5}{8}]\}}\lambda\quad\text{and}\quad}
\tau\yeq\tau(H,\theta)
\yeq 
{\sred1_{\{H\in(\half,\frac{5}{8}]\}}}{\bf G}(\theta)^{-1}b_\infty(\theta)-\beta(\theta). 
\eea
{\bred Recall that the constant $\lambda$ is defined in (\ref{202308150236}).}
More precisely, 
\bea\label{202308141933}
g(x) 
&=&
\phi(x;0,c_0)\big\{
-2\kappa x
+(\kappa x^2+\tau)H_1(x,c_0)\big\}
\nn\\&=&
\phi(x;0,c_0)\big\{
(\tau-2\kappa c_0)H_1(x,c_0)+\kappa x^2H_1(x,c_0)\big\}
\nn\\&=&
\phi(x;0,c_0)\big\{
(\tau-2\kappa c_0)H_1(x,c_0)+\kappa c_0^2H_3(x,c_0)+3\kappa c_0H_1(x,c_0)\big\}
\nn\\&=&
\phi(x;0,c_0)\big\{
(\tau+\kappa c_0)H_1(x,c_0)+\kappa c_0^2H_3(x,c_0)\big\}
\eea
Remark that $H_3(x,c_0)=c_0^{-3}x^3-3c_0^{-2}x$ and 
\beas 
x^2H_1(x,c_0) 
&=& 
c_0^2H_3(x;0,c_0)+3c_0H_1(x;0,c_0).
\eeas
With 
$\tau$ and $\kappa$ of (\ref{202308150237}) 
and 
$c_3'$ of (\ref{202308141444}), 
set 
\bea\label{202308141938}
c_1 \yeq \tau+\kappa c_0 
\quad\text{and}\quad
c_3\yeq c_3'+3{\bred\lambda} c_0^2.
\eea
{\bred Remark that 
$1_{\{H\in(\half,\frac{5}{8}]\}}c_3\yeq 1_{\{H\in(\half,\frac{5}{8}]\}}(c_3'+3\kappa c_0^2)$.} 
Then the resulting asymptotic expansion formula for $F_T^{2C_1}$ is given by 
$p_{H,T}$ of (\ref{202308141943}). 

Since the estimator $\wh{\theta}_T$ takes values in the bounded set $\Theta$ 
and as already mentioned $\psi_T^{2C_1}-1=O_{D_\infty}(T^{-L})$ for every $L>0$, 
it is easy to show 
\beas 
\sup_{g\in\cale(a,b)}\big|E\big[g\big(T^{1/2}(\wh{\theta}_T-\theta)\big)\big]
-
E\big[g\big(F^{2C_1}_T\big)\big]\big|
&=&
O(T^{-L})\quad(T\to\infty)
\eeas
for every $L>0$. 
Thus, we obtain the asymptotic expansion and its error bound for $T^{1/2}(\wh{\theta}_T-\theta)$. 
\qed\halflineskip

{\bred 
\subsection{Proof of Theorem \ref{202403210212}}
Define $c_{1,1}^+$ and $c_{1,2}^+$ as 
\bea\label{202403210242}
c_{1,1}^+\yeq{\bf G}(\theta)^{-1}b_\infty(\theta)+\lambda c_0
\quad\text{and}\quad
c_{1,2}^+\yeq-\beta(\theta).
\eea
Then, by the definition (\ref{202403210310}) of $\p_{H,T}^+$ and the argument in Section \ref{202403210312}, 
we see 
\beas
\sup_{g\in\cale(a,b)}\bigg|\int_\bbR g(x)\big(p_{H,T}(x)-p_{H,T}^+(x)\big)dx\bigg|
&=&
o(T^{-\sfq(H)})
\eeas
as $T\to\infty$, for every $a,b>0$. 
Therefore, 
(\ref{202403210214}) follows from (\ref{202403210322}) of Theorem \ref{202308141824}. 
\qed\halflineskip
}

\begin{en-text}
\section{Appendix}
We show $\sfB_m(2H-2,...,2H-2)<\infty$ when $H>\frac{m+1}{2m}$. 
Recall the formula
\beas 
\int_x^z(z-y)^p(y-x)^qdy &=& (z-x)^{p+q{\colorr+1}}\sfB(p+1,q+1)\qquad(x<z;\>p,q>-1).
\eeas
We have $H_0:=2H-2>-1$ since $H>1/2$, and 
\beas
&&
\int_{\{0<x_m<x_{m-1}<\cdots<x_1<1\}}|x_1-x_2|^{H_0}|x_2-x_3|^{H_0}\cdots
|x_{m-1}-x_m|^{H_0}|x_m-x_1|^{H_0}
dx_1...dx_m
\\&=&
\sfB(H_0+1,H_0+1)
\int_{\{0<x_{m-1}<\cdots<x_1<1\}}|x_1-x_2|^{H_0}|x_2-x_3|^{H_0}\cdots
|x_{m-1}-x_1|^{2H_0+1}
dx_1...dx_{m-1}
\\&=&
\sfB(H_0+1,H_0+1)\sfB(H_0+1,2H_0+2)
\int_{\{0<x_{m-2}<\cdots<x_1<1\}}|x_1-x_2|^{2H-2}|x_2-x_3|^{2H-2}\cdots
\nn\\&&\hspace{10pt}
\cdots|x_{m-2}-x_1|^{3H_0+2}dx_1...dx_{m-2}
\\&=&
\sfB(H_0+1,H_0+1)\sfB(H_0+1,2H_0+2)
\int_{\{0<x_{m-2}<\cdots<x_1<1\}}|x_1-x_2|^{2H-2}|x_2-x_3|^{2H-2}\cdots
\nn\\&&\hspace{10pt}
\cdots|x_{m-2}-x_1|^{3H_0+2}dx_1...dx_{m-2}
\\&=&\cdots
%
\\&=&
\sfB(H_0+1,H_0+1)\sfB(H_0+1,2H_0+2)\cdots \sfB(H_0+1,(m-2)H_0+(m-2))
\nn\\&&\times
\int_{\{0<x_2<x_1<1\}}|x_1-x_2|^{H_0}|x_2-x_1|^{(m-1)H_0+(m-2)}
dx_1dx_2
\\&=&
\sfB(H_0+1,H_0+1)\sfB(H_0+1,2(H_0+1))\cdots \sfB(H_0+1,(m-2)(H_0+1))
\nn\\&&\times
\int_{\{0<x_2<x_1<1\}}|x_2-x_1|^{mH_0+m-2}dx_1dx_2
\\&=&
\sfB(H_0+1,H_0+1)\sfB(H_0+1,2(H_0+1))\cdots \sfB(H_0+1,(m-2)(H_0+1))
\nn\\&&\times
(mH_0+m-1)^{-1}
\int_{\{0<x_1<1\}}x_1^{mH_0+m-1}dx_1dx_2
\quad\underline{\text{if }mH_0+m-2>-1}
\\&=&
\sfB(H_0+1,H_0+1)\sfB(H_0+1,2(H_0+1))\cdots \sfB(H_0+1,(m-2)(H_0+1))
\nn\\&&\times
(mH_0+m-1)^{-1}(mH_0+m)^{-1}\><\infty
\eeas
if $mH_0+m-1>0$, or $H>\frac{m+1}{2m}$. 
\end{en-text}
\begin{en-text}
In particular, for $m=2$, if $H>3/4$ (not used however), then 
\beas 
\lim_{T\to\infty}T^{3-4H}k_2(F_T) 
&=& 
8\{H^2(2H-1)^2/2\}\sfB(2H-1,2H-1)(4H-3)^{-1}(4H-2)^{-1}H^2(2H-1)^2/2
\nn\\&=&
\sfB(2H-1,2H-1)(4H-3)^{-1}H^2(4H-2)^2. 
\eeas
For $m=3$, if $H>2/3$, then 
\beas
\lim_{T\to\infty}T^{\frac{9}{2}-6H}k_3(F_T)
&=& 
8\cdot H^3(2H-1)^3\cdot 6\cdot \sfB(2H-1,2H-1)\sfB(2H-1,2(2H-1))(6H-4)^{-1}(6H-3)^{-1}
\nn\\&=&
12H^3(2H-1)(3H-2)^{-1}\sfB(2H-1,2H-1)\sfB(2H-1,2(2H-1)).
\eeas
For $m=4$, if $H>5/8$, then 
\beas
\lim_{T\to\infty}T^{6-8H}k_4(F_T)
&=& 
16\cdot 3H^4(2H-1)^4\cdot 24\cdot \sfB(2H-1,2H-1)\sfB(2H-1,2(2H-1))\sfB(2H-1,3(2H-1))
\nn\\&&\times
(8H-5)^{-1}(8H-4)^{-1}
\nn\\&=&
288H^4(2H-1)^3(8H-5)^{-1}\sfB(2H-1,2H-1)\sfB(2H-1,2(2H-1))\sfB(2H-1,3(2H-1))
\nn\\&=&
288H^4(2H-1)^3(8H-5)^{-1}\{\Gamma(2H-1)\}^4\{\Gamma(4(2H-1))\}^{-1}. 
\eeas
\end{en-text}

{\bred
\section{Simulation study}\label{202403280137}
The performance of the asymptotic expansion formula $p_{H,T}$ of (\ref{202308141943}) will be investigated 
by simulations. 
We consider the parameter values $\theta=2$ and $H\in\{0.55,0.625,0.7\}$. 
The number of replications in each Monte Carlo simulation is $10^5$. 
The YUIMA package (cf. \cite{Yuima2014,iacus2017simulation}) is used for the study.

Figure \ref{202403280353} shows the asymptotic expansion formula $p_{0.55,50}$  
captures the skewness of the distribution of the estimation error in the time horizon $T=50$. 
On the other hand, the normal approximation improves for $T=100$ as in Figure \ref{202403280354}. 
\begin{figure}[H]
\hspace{0.04\columnwidth}
  \begin{minipage}[b]{0.40\columnwidth}
    \centering
    \includegraphics[clip,width=\columnwidth]{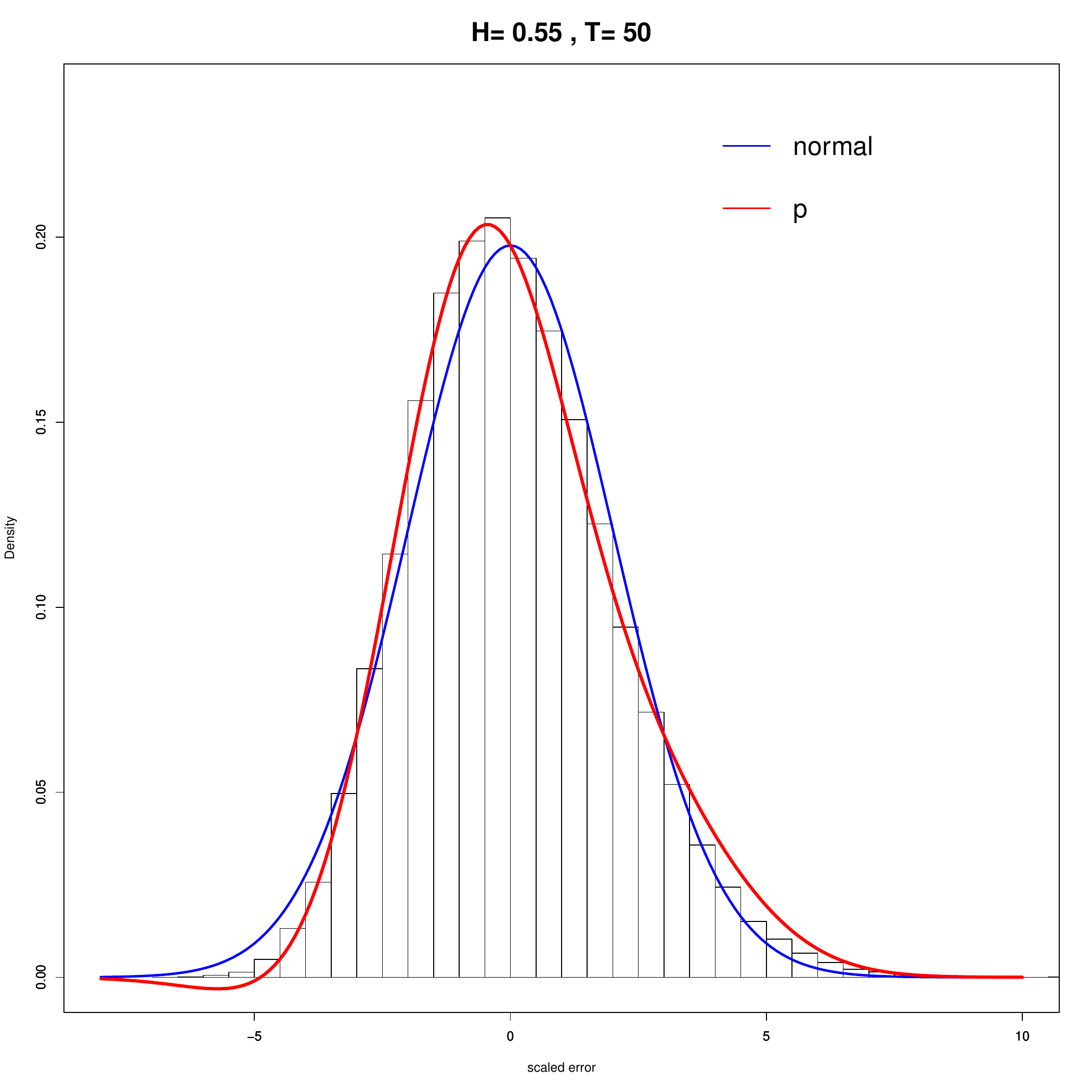}
    \caption{$N(0,c_0)$ and $p_{0.55,50}$}
     \label{202403280353}
  \end{minipage}
  \hspace{0.04\columnwidth} 
  \begin{minipage}[b]{0.40\columnwidth}
    \centering
    \includegraphics[clip, width=\columnwidth]{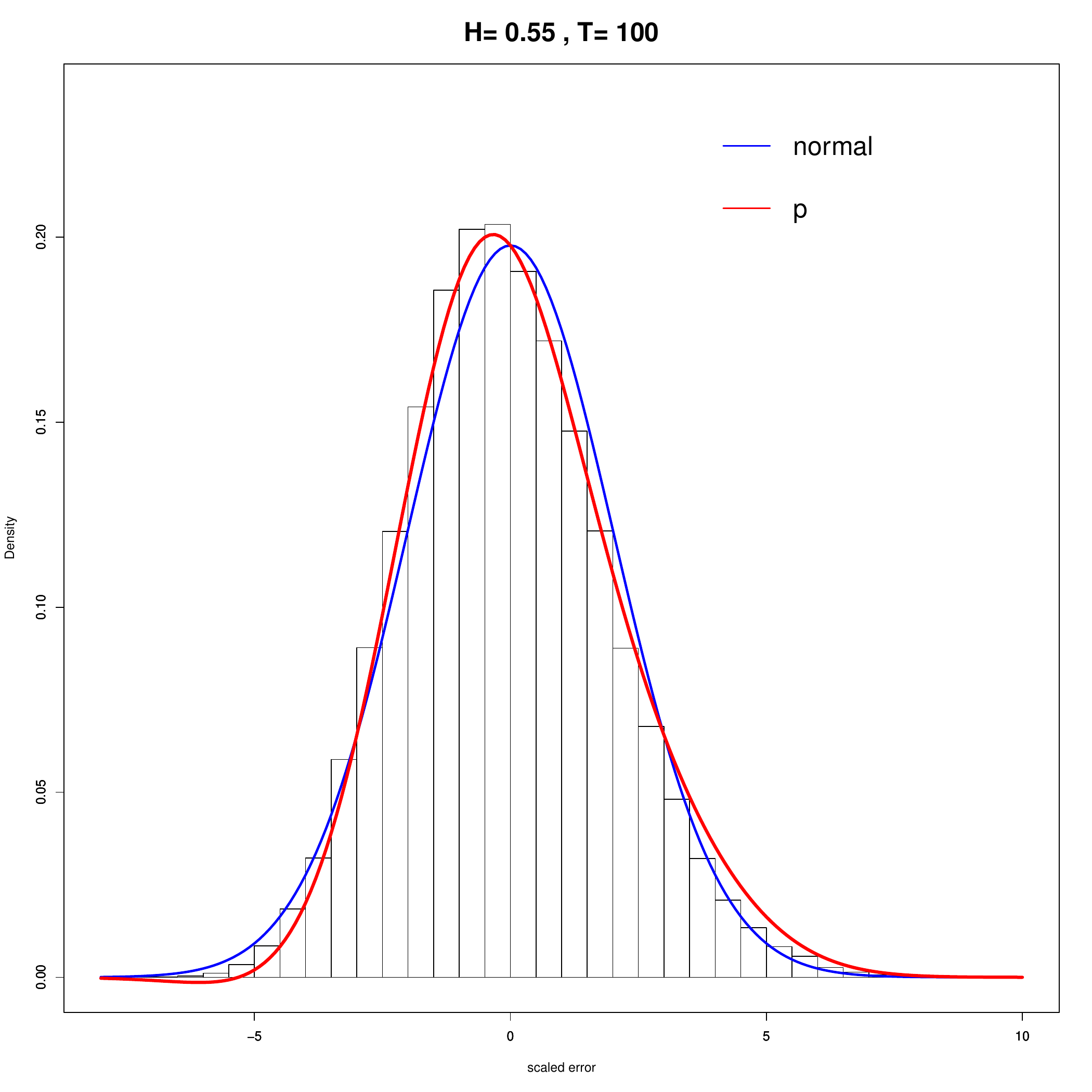}
    \caption{$N(0,c_0)$ and $p_{0.55,100}$}
    \label{202403280354}
  \end{minipage}
  \hspace{0.08\columnwidth}
\end{figure}

The value $H=5/8=0.625$ is the threshold of $T$'s exponents $-1/2$ and $4H-3$ of the first-order correction term 
of the asymptotic expansion. 
Figures \ref{202403280255} and \ref{202403280256} show that 
the asymptotic expansion formulas have caught the skewness of the distribution. 
The correction becomes smaller for the larger $T$. 
Since the first-order correction by the asymptotic expansion consists of the two terms, 
it is a bit unexpected that the difference between the histogram and the normal distribution is rather small. 
However, it is natural in a sense because the relative effect of the skewness decreases down toward $5/8$ on $(1/2,5/8]$,  
and the relative effect of the gap between the real variance and $c_0$ goes down toward $5/8$ on $[5/8,3/4)$. 
\begin{figure}[H]
\hspace{0.04\columnwidth}
  \begin{minipage}[b]{0.40\columnwidth}
    \centering
    \includegraphics[clip,width=\columnwidth]{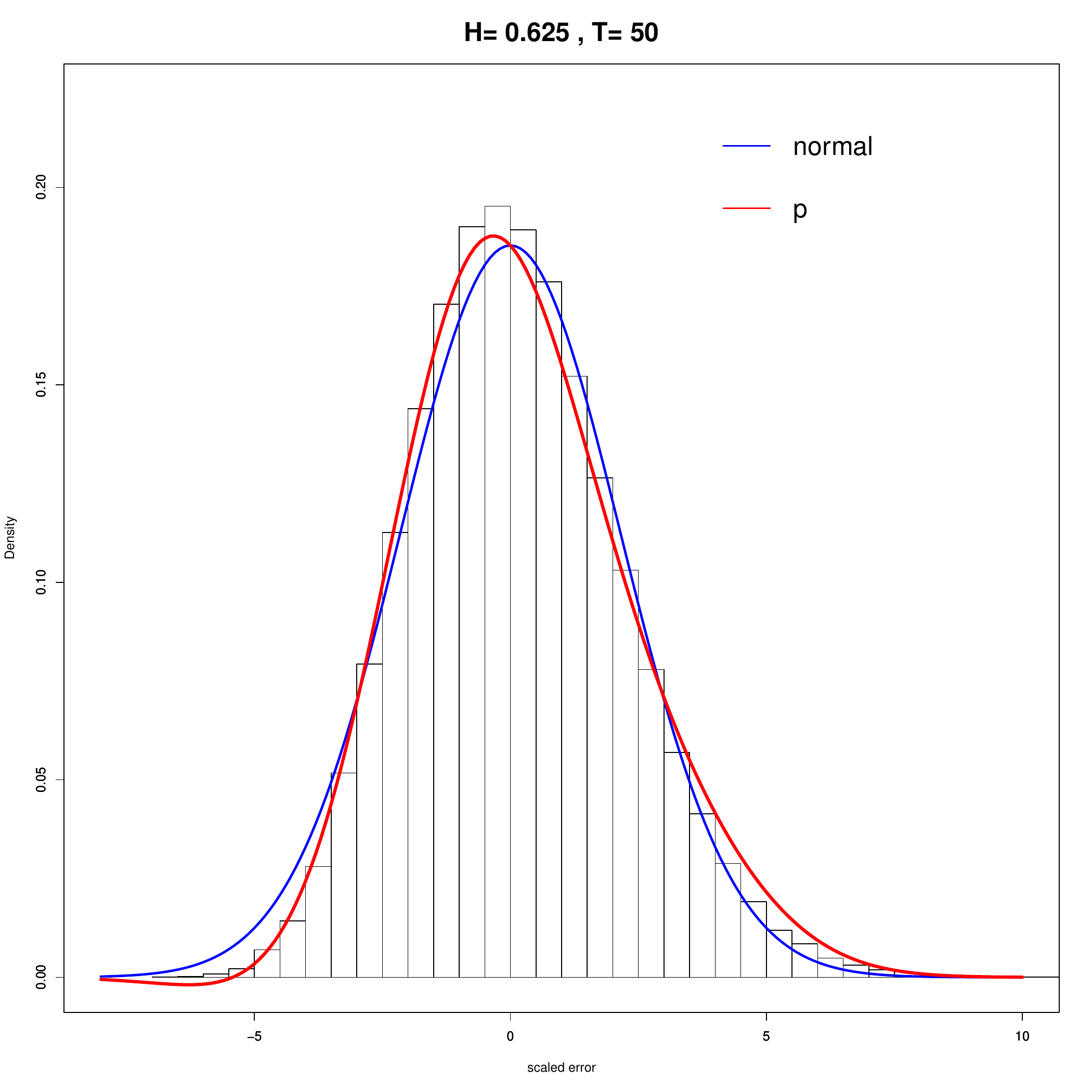}
    \caption{$N(0,c_0)$ and $p_{0.625,50}$}
     \label{202403280255}
  \end{minipage}
  \hspace{0.04\columnwidth} 
  \begin{minipage}[b]{0.40\columnwidth}
    \centering
    \includegraphics[clip, width=\columnwidth]{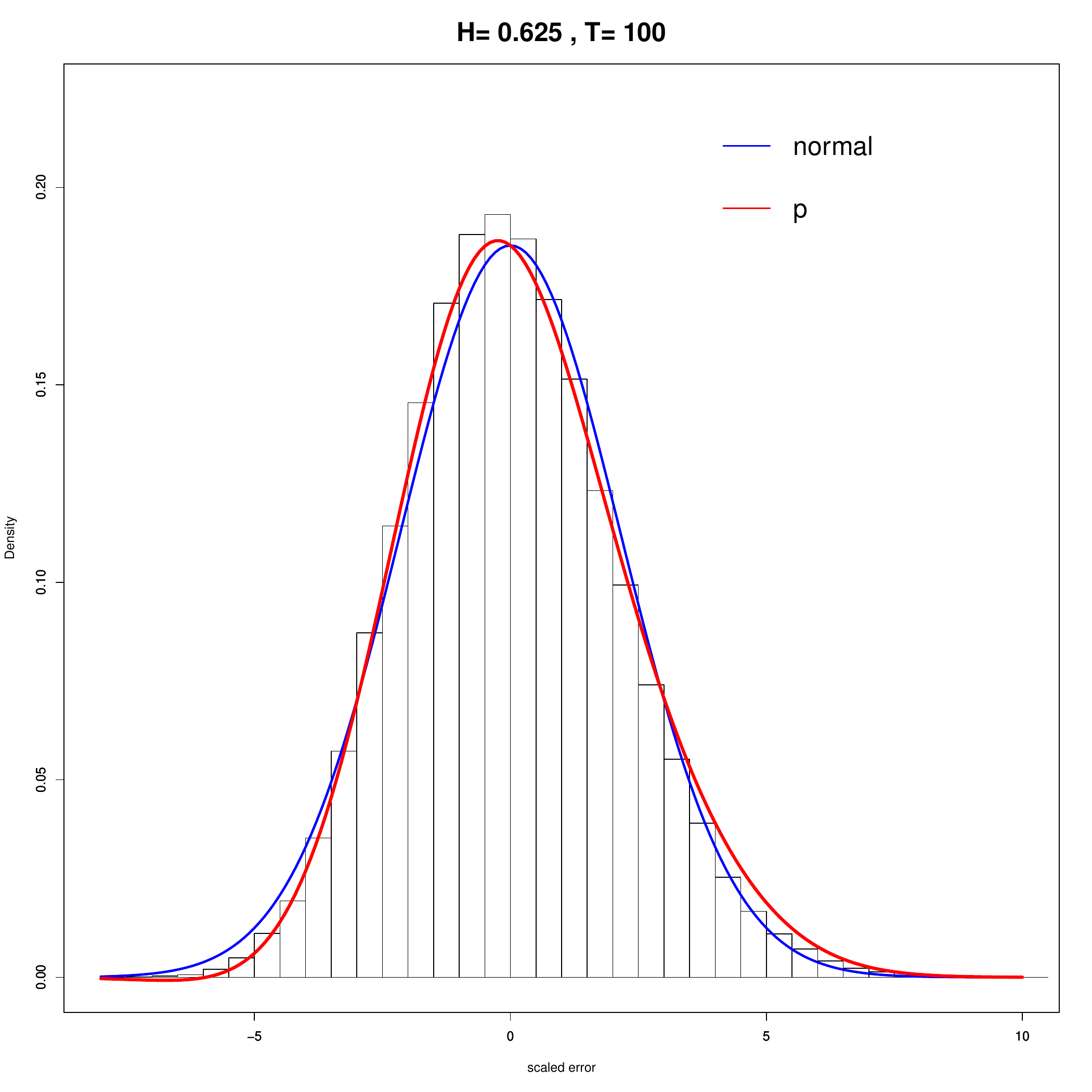}
    \caption{$N(0,c_0)$ and $p_{0.625,100}$}
    \label{202403280256}
  \end{minipage}
  \hspace{0.08\columnwidth}
\end{figure}

In the case $H=0.7$, 
Figure \ref{202403280257} shows the asymptotic expansion fairly improves the normal approximation. 
However, some discrepancy remains yet between the asymptotic expansion and the histogram, even for $T=100$, 
for which the normal approximations performed better when $H=0.55$ and $0.625$, as observed above. 
The value $H=0.7$ is near to the upper bound of the interval $(1/2,3/4)$ (more generally (0,3/4)) 
of $H$ for the valid normal approximation with the scaling $T^{1/2}$. 
Hu et al. \cite{hu2019parameter} showed that 
the limit becomes a normal distribution for $H=3/4$ with the rate of convergence $T^{1/2}/\sqrt{\log T}$, 
and a Rosenblatt distribution if $H$ exceeds $3/4$ with the rate $T^{2-2H}$. 
This fact explains the relatively large discrepancy between the histogram and the normal approximation under rate $T^{1/2}$. 
The asymptotic expansion is trying to approximate the histogram, while it still has a gap 
since the first-order asymptotic expansion $p_{0.7,100}$ does not incorporate the effect of the kurtosis nor the higher-order moments of the variable. 
The approximations by the asymptotic expansion and normal distribution are improved when $T=400$ 
as Figure \ref{202403280258} though the error of the normal approximation is not small yet. 
\begin{figure}[H]
\hspace{0.04\columnwidth}
  \begin{minipage}[b]{0.40\columnwidth}
    \centering
    \includegraphics[clip,width=\columnwidth]{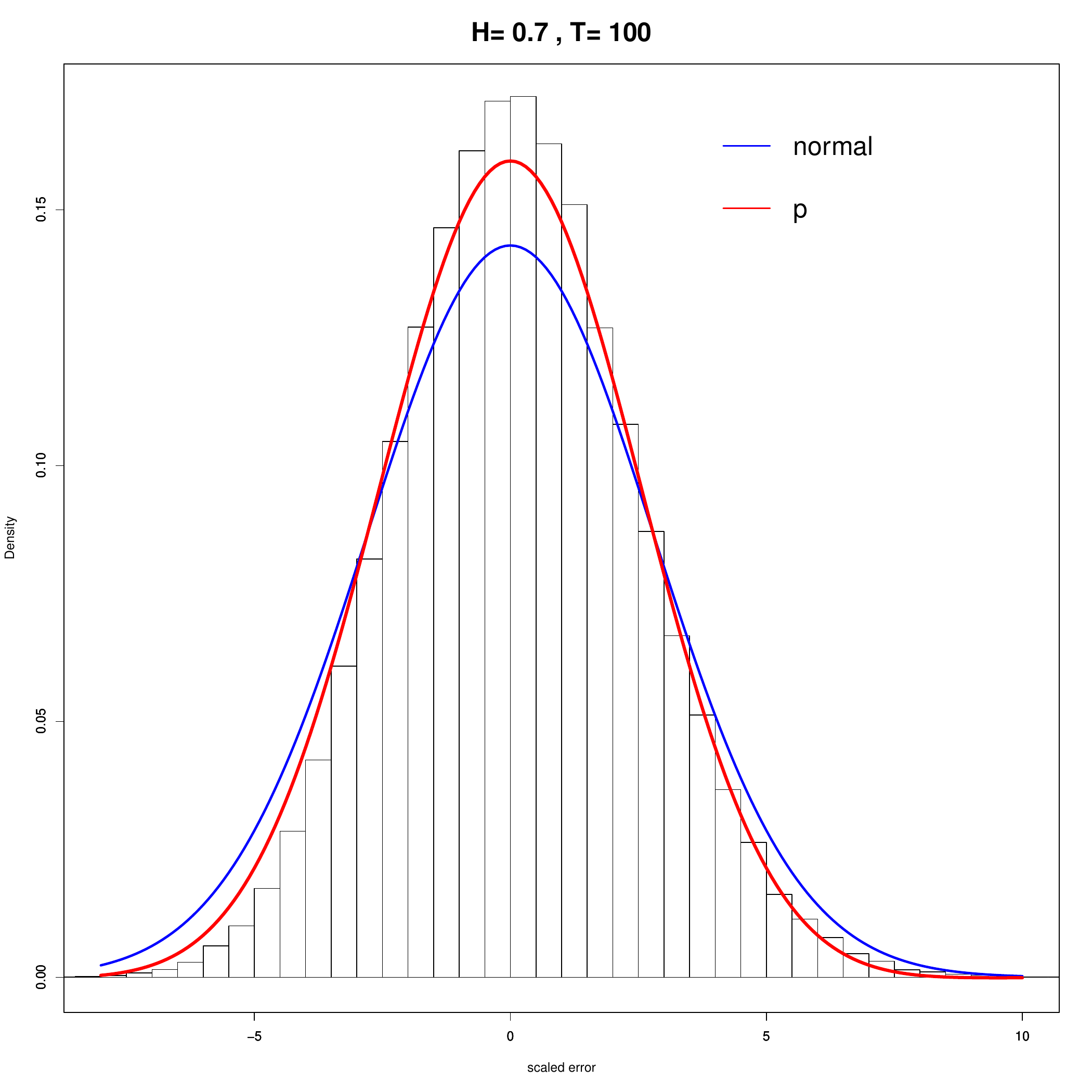}
    \caption{$N(0,c_0)$ and $p_{0.7,100}$}
     \label{202403280257}
  \end{minipage}
  \hspace{0.04\columnwidth} 
  \begin{minipage}[b]{0.40\columnwidth}
    \centering
    \includegraphics[clip, width=\columnwidth]{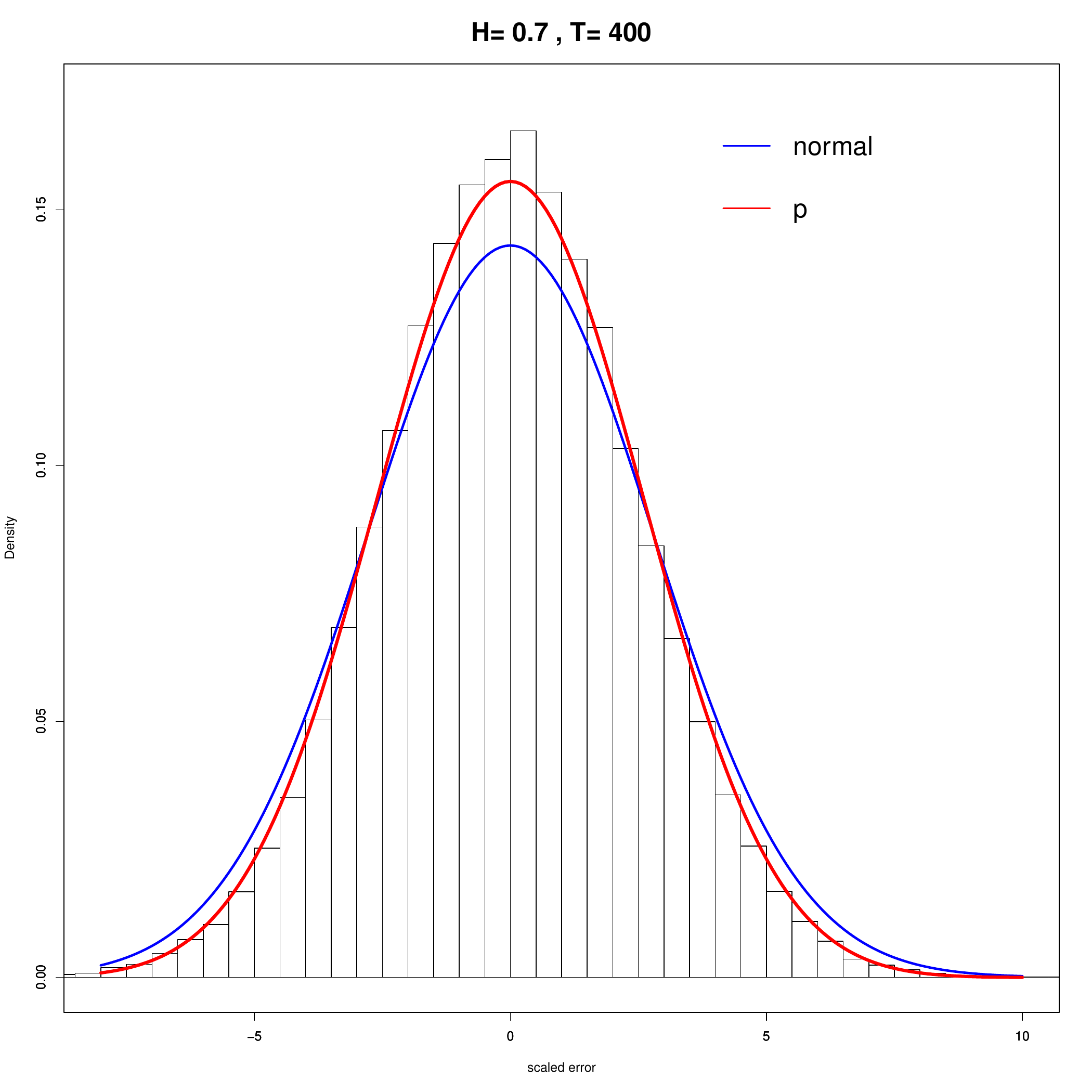}
    \caption{$N(0,c_0)$ and $p_{0.7,400}$}
    \label{202403280258}
  \end{minipage}
  \hspace{0.08\columnwidth}
\end{figure}

}

\bibliographystyle{spmpsci}      
\bibliography{bibtex-20210212+-20240210+}   

\end{document}